\documentclass[12pt,a4paper]{article}
\usepackage{amsfonts}
\usepackage{amsmath,amssymb,amsthm}
\usepackage{amssymb,latexsym}

\newtheorem{theorem}{Theorem}[section]
\newtheorem{lemma}{Lemma}[section]
\newtheorem{prop}{Proposition}[section]
\newtheorem{cor}{Corollary}[section]

\newcommand{\LP}{L\'{e}vy process}
\newcommand{\R}{\mathbb{R}}
\newcommand{\Rd}{\mathbb{R}^{d}}
\newcommand{\nN}{n \in \mathbb{N}}

\newcommand{\C}{\mathbb{C}}
\newcommand{\E}{\mathbb{E}}

\newcommand{\ds}{\displaystyle}
\newcommand{\g}{\mathfrak{g}}

\newcommand{\Dom}{\mbox{Dom}}

\newcommand{\bean}{\begin{eqnarray*}}
\newcommand{\eean}{\end{eqnarray*}}

\newcommand{\la}{\langle}
\newcommand{\ra}{\rangle}

\newcommand{\Z}{\mathbb{Z}}

\newcommand{\Ad}{\mbox{Ad}}
\newcommand{\at}{\mathfrak{a}^{*}}
\newcommand{\Ran}{\mbox{Ran}}

\newcommand{\fk}{\mathfrak{k}}
\newcommand{\fp}{\mathfrak{p}}
\newcommand{\fa}{\mathfrak{a}}
\newcommand{\fn}{\mathfrak{n}}

\newcommand{\Exp}{\mbox{Exp}}

\date{}

\title {A Generalised Gangolli-L\'{e}vy-Khintchine Formula for Infinitely Divisible Measures and
L\'{e}vy Processes on Semi-Simple Lie Groups and Symmetric Spaces}

\author{David Applebaum\footnote{D.Applebaum@sheffield.ac.uk}\\
 School of Mathematics and Statistics,\\ University of Sheffield, \\
Sheffield S3 7RH\\
United Kingdom.\\
 \and
    Anthony Dooley\footnote{a.h.dooley@bath.ac.uk}{}\\
 Department of Mathematical Sciences,\\
University of Bath,\\Bath BA2 7AY\\
United Kingdom.}

\begin{document}
\maketitle

\begin{abstract} In 1964 R.Gangolli published a L\'{e}vy-Khintchine type formula which characterised $K$ bi-invariant infinitely divisible probability measures on a symmetric space $G/K$. His main tool was Harish-Chandra's spherical functions which he used to construct a generalisation of the Fourier transform of a measure. In this paper we use generalised spherical functions (or Eisenstein integrals) and extensions of these which we construct using representation theory  to obtain such a characterisation for arbitrary infinitely divisible probability measures on a non-compact symmetric space. We consider the example of hyperbolic space in some detail.
\end{abstract}

\section{Introduction} \label{sec1}

A probability measure on a topological group $G$ is said to be {\it infinitely divisible} if it has a convolution $n$th root for all natural numbers $n$. A stochastic process taking values in $G$ is called a {\it L\'{e}vy process} if it has stationary and independent increments and is stochastically continuous. The probability laws of the process then form a (weakly continuous) convolution semigroup of measures. If $G$ is a connected Lie group having at least one finite dimensional linear representation that has a discrete kernel then the {\it Dani-McCrudden embedding theorem} states that every infinitely divisible measure is of the form $\mu_{1}$ for some convolution semigroup of measures $(\mu_{t}, t \geq 0)$ (see \cite{DM, DM1}.)

The investigation of convolution semigroups (and hence L\'{e}vy processes) on Lie groups dates back to Hunt in 1956 who characterised them via the infinitesimal generators of the associated semigroup of operators acting on a Banach space of continuous functions on $G$ \cite{Hu}. Since then there has been much progress in developing understanding of these objects and the reader is directed to the monograph \cite{Liao} for insight. In the case where $G = \Rd$, a more direct characterisation of infinitely divisible measures $\mu$ is obtained by using the Fourier transform to derive the classical {\it L\'{e}vy-Khintchine formula}
\begin{equation} \label{LK1}
\int_{\Rd}e^{i u \cdot x}\mu(dx) = e^{-\eta(u)},
\end{equation}
where $\eta: \Rd \rightarrow \C$ is a continuous, hermitian negative-definite function (see e.g. \cite{Sa, Abook} for details.)

Since any globally Riemannian symmetric space $M$ is a homogeneous space $G/K$ where $G$ is a connected semisimple Lie group and $K$ is a compact subgroup, L\'{e}vy processes and infinitely divisible measures on $M$ can be defined to be the images of those on $G$ under the associated canonical surjection. Hence probability theory on $M$ is reduced to that on the Riemannian symmetric pair $(G,K)$. In 1964 Gangolli \cite{Gang1} found a precise analogue of the L\'{e}vy-Khintchine formula (\ref{LK1}) in this context for $K$ bi-invariant infinitely divisible measures $\mu$ on $G$. A key ingredient was the generalisation of the Fourier transform in (\ref{LK1}) to the spherical transform formed using a spherical function $\Phi$ on $G$ (see e.g. \cite{Hel2} for background on these.) In fact we then have
\begin{equation} \label{LK2}
\int_{G}\Phi(g)\mu(dg) = e^{-\eta_{\Phi}},
\end{equation}
where $\eta_{\Phi} \in \C$. If $G$ is semisimple and $M$ is irreducible then
$$ \eta_{\Phi} = a \lambda_{\Phi} + \int_{G - \{e\}}(1 - \Phi(g))\nu(dg),$$
where $\lambda_{\phi}$ is an eigenvalue of the Laplace-Beltrami operator on $G$ (corresponding to the eigenvector $\Phi$), $a \geq 0$ and $\nu$ is a L\'{e}vy measure on $G - \{e\}$ (see also \cite{App0} and \cite{LW}.)

Gangolli's paper \cite{Gang1} stimulated a great deal of further work on this subject. In \cite{Gang2} he investigated the sample paths of L\'{e}vy processes on $M$ (see also \cite{App0}.) Integrability and transience have been studied in \cite{Berg}, \cite{BeFa}, \cite{Hey2} and \cite{App5}. Subordination on symmetric spaces was investigated in \cite{App4} and \cite{AlGo}  while conditions for the existence of smooth densities were found in \cite{LW}. This work also inspired generalisations to general Riemannian manifolds \cite{AE}, Gelfand pairs \cite{Hey1} and hypergroups \cite{BlHe}.

 We are particularly interested in the non-compact case. Then the spherical functions on $G$ take the form
 \begin{equation} \label{sphernc}
\Phi_{\lambda}(\sigma) = \int_{K}e^{(i\lambda +
\rho)(A(k\sigma))}dk \end{equation}
for all $\sigma \in G$, where $A$ is the abelian part of the Iwasawa decomposition of $G$, $\rho$ is the half-sum of positive roots and the ``parameter'' $\lambda$ runs through the (real) dual space of the abelian part of the Iwasawa decomposition of the Lie algebra of $G$ (see below for more details.) Our goal in this paper is to extend Gangolli's result to general infinitely divisible measures on $G$ (without any bi-invariance assumptions.)  Spherical functions are no longer adequate tools for us to build the Fourier transform that we require. We note that these functions have been generalised to Eisenstein integrals (or generalised spherical functions) $\Phi_{\lambda, \pi}$ \cite{Hel3} wherein the measure $dk$ in (\ref{sphernc}) is replaced by the ``matrix-valued weight'' $\pi(k)dk$ where $\pi$ is an irreducible representation of the compact group $K$. In fact we find that this generalisation is not enough. To derive our formula we need to take a further step and consider a double parametrisation by the unitary dual $\widehat{K}$ to form objects $(\Phi_{\lambda, \pi, \pi^{\prime}}, \pi, \pi^{\prime} \in \widehat{K})$ which we regard as infinite matrix valued functions defined on $G$. The construction of these {\it generalised Eisenstein integrals} utilises representation theory techniques and is carried out in section 3 of this paper. This part of the work may be of independent interest to specialists in harmonic analysis on Lie groups and symmetric spaces.

In section 5 of the paper we derive our L\'{e}vy-Khintchine type formula for semisimple groups. In fact we show that if $\mu$ is an infinitely divisible probability measure on $G$ then
\begin{equation} \label{LK3}
\int_{G}\Phi_{\lambda}(g)\mu(dg) = \Exp(\psi_{\lambda}),
\end{equation}
and we find a precise form for the infinite matrix $(\psi_{\lambda, \pi, \pi^{\prime}}, \pi, \pi^{\prime} \in \widehat{K})$. The meaning of the ``exponential'' Exp for an infinite matrix is made precise in section 4.

Our approach is different from Gangolli's. He worked directly with infinitely divisible measures. We work with convolution semigroups (in a similar way to \cite{LW} and \cite{App0}) so that we can take advantage of Hunt's classification \cite{Hu} and first obtain (\ref{LK3}) within that context. We then obtain the result for infinitely divisible measures from the Dani-McCrudden embedding theorem as described above. In section 6 we apply this theory to symmetric spaces $G/K$, and in section 7 we give a detailed description of the results on hyperbolic space, relating it to the Helgason Fourier transform. As pointed out to us by Peter Kim, these latter results are applicable to the statistical problem of deconvolution density estimation in hyperbolic space where they enable the assumption of spherical symmetry of the error density in \cite{HKKM} to be dropped. More generally, the results in this paper have potential for application to the study of transience and recurrence of L\'{e}vy processes in symmetric spaces (see e.g. \cite{App5}) and to the study of limiting properties and rates of convergence of components of L\'{e}vy processes in groups and symmetric spaces (see Chapters 6 and 7 of \cite{Liao}.)

As we are aiming this paper at an audience of both probabilists and harmonic analysts, we have included, in section 8 an appendix where we treat the properties of induced representations, principal series of semisimple groups and their $K$-types, which are needed as background for the present work. Whilst this material is in some senses known, it is hard to access in a uniform and consistent presentation. We hope that it will engender a wider readership of this work.

\vspace{5pt}

{\it Notation.} If $G$ is a Lie group then ${\cal B}(G)$ is its Borel $\sigma$-algebra and $C_{0}(G)$ is the Banach space of real-valued continuous functions on $G$ which vanish at infinity, equipped with the usual supremum norm. If $T$ is a linear operator acting on $C_{0}(G)$ then Dom$(T)$ denotes its maximal domain. The $*$-algebra of all bounded linear operators on a complex separable Hilbert space $H$ is denoted ${\cal L}(H)$. Einstein summation convention is used throughout this paper. If $f:G \rightarrow \C$ is continuous we write $f^{\#}(g) = f(g^{-1})$ for $g \in G$. The (algebraic) dual of a complex vector space $V$ will be denoted $V^{*}$.

\section{L\'{e}vy Processes on Lie Groups}\label{sec2}

Let $G$ be a Lie group with neutral element $e$, Lie algebra $\g$ and dimension $n$ and let $Z = (Z(t), t \geq 0)$ be a $G$-valued stochastic process defined on some probability space $(\Omega, {\cal F}, P)$. The right increment of $Z$ between the times $s$ and $t$ where $s \leq t$ is the random variable $Z(s)^{-1}Z(t)$. We say that $Z$ is a \LP~on $G$ if $Z$ has stationary and independent right-increments, is stochastically continuous and $Z(0) = e$ (a.s.). Each $Z(t)$ has law $\mu_{t}$ and $(\mu_{t}, t \geq 0)$ is a weakly continuous convolution semigroup of probability measures with $\mu_{0}$ being Dirac mass at $e$. Conversely given any such convolution semigroup we can always construct a L\'{e}vy process $Z$ for which the law of $Z(t) $ is $\mu_{t}$ for each $t \geq 0$ on the space of all paths from $[0, \infty)$ to $G$ by using Kolmogorov's construction theorem (see e.g. \cite{App5} and Theorem 10.5 on pp.55-7 of \cite{Sa} for a detailed account of the case $G=\R^{n}$.) We refer the reader to \cite{Liao} for a monograph treatment of Lie group valued L\'{e}vy processes.  Given an arbitrary \LP~$Z$ we define an associated $C_{0}$-contraction semigroup $(T_{t}, t \geq 0)$ on $C_{0}(G)$ by the prescription:

\begin{equation} \label{sg}
(T_{t}f)(g) = \E(f(gZ(t))) =
\int_{G}f(gh)\mu_{t}(dh),
\end{equation}
for each $f \in C_{0}(G), g \in G, t \geq 0$.

The infinitesimal generator ${\cal L}$  of this semigroup was
characterised by Hunt \cite{Hu} in 1956 (see \cite{Liao} pp.52-61 for a more recent treatment).  We first fix a basis $(X_{j}, 1
\leq j \leq n)$ of $\g$ and define the dense linear manifold
$C^{2}(G)$ by
$$ C^{2}(G) = \{ f \in C_{0}(G); X_{i}(f) \in C_{0}(G)~ \mbox{and}~
X_{i}X_{j}(f) \in C_{0}(G) ~\mbox{for all}~ 1 \leq i,j \leq n\},$$
where the Lie algebra acts as left-invariant vector fields.

There exist functions ${x_{i} \in C_{c}^{\infty}(G), 1 \leq i \leq
n}$ so that $(x_{1}, \ldots, x_{n})$ are a system of canonical
co-ordinates for $G$ at $e$. A measure $\nu$ defined on ${\cal
B}(G - \{e\})$ is called a {\it L\'{e}vy measure} whenever
     $$  \int_{G-\{e\}}\left(\sum_{i=1}^{n}x_{i}(g)^{2}\right) \nu(dg) < \infty~\mbox{and}~\nu(U^{c}) < \infty,$$ for any neighbourhood $U$ of $e$.

\begin{theorem}[Hunt's theorem] \label{Hunt}
Let $(\mu_{t}, t \geq 0)$ be a weakly continuous convolution semigroup of measures
in $G$ with infinitesimal generator ${\cal L}$ then
\begin{enumerate}
\item $C^{2}(G) \subseteq \mbox{Dom}({\cal L})$. \item For each
$g \in G, f \in C^{2}(G)$,
\begin{equation} \label{hu}
{\cal L}f(g)  = b_{i}X_{j}f(g) +
a^{ij}X_{i}X_{j}f(g)
   + \int_{G-\{e\}}(f(gh) - f(g) -
   x^{i}(h)X_{i}f(g))\nu(dh),
\end{equation}
where $b = (b^{1}, \ldots b^{n}) \in {\R}^{n}, a = (a^{ij})$ is a
non-negative-definite, symmetric $n \times n$ real-valued matrix
and $\nu$ is a L\'{e}vy measure on $G-\{e\}$.

\end{enumerate}
Conversely, any linear operator with a representation as in
(\ref{hu}) is the restriction to $C^{2}(G)$ of the infinitesimal generator of a semigroup of convolution operators
on $C_{0}(G)$ that is induced by a unique weakly continuous convolution semigroup of probability measures.
\end{theorem}

From now on we assume that the Lie group $G$ is connected, semi-simple and has a finite centre. Let $K$ be a maximal compact subgroup of $G$ with Lie algebra $\fk$ and note that $(G,K)$ is a non-compact Riemannian symmetric pair. The Iwasawa decomposition gives a global diffeomorphism
between $G$ and a direct product $NAK$ where $A$ and $N$ are
simply connected with $A$ being abelian and $N$ nilpotent wherein
each $g \in G$ is mapped onto $n(g)\exp(A(g))u(g)$ where $u(g) \in
K, n(g) \in N$ and $A(g) \in \fa$ which is the Lie algebra
of $A$. Let $\widehat{K}$ be the unitary dual of $K$, i.e. the set of all equivalence classes (modulo unitary transformations) of irreducible
unitary representations of $K$. As is standard, we will frequently identify classes in $\widehat{K}$ with representative elements.  We recall that the {\it generalised
spherical functions} or {\it Eisenstein integrals} \cite{Hel3} are the
matrix-valued functions:
\begin{equation} \label{gsph}
\Phi_{\lambda, \pi}(g) = \int_{K}e^{(i\lambda +
\rho)(A(kg))}\overline{\pi}(k)dk \end{equation}
for $g \in G$, where $\lambda \in \fa^{*}_{\C}$ which
is the complexification of the dual space $\fa^{*}$ of
$\fa$ and $\rho$ is half the sum of positive roots
(relative to a fixed lexicographic ordering.) Note that $\overline{\pi}$ is the conjugate of $\pi \in \widehat{K}$ which we will discuss in detail within the next section.

One approach to generalising Gangolli's L\'{e}vy-Khintchine formula \cite{Gang1} (see also \cite{App0}, \cite{LW})
would be to seek to find the (principal part of the) logarithm of
\begin{equation} \label{log1}
\alpha_{\lambda, \pi}(t) := T_{t}(\Phi_{\lambda, \pi})(e) =
\int_{G}\Phi_{\lambda, \pi}(g)\mu_{t}(dg).
\end{equation}
 We will see that this plan is not adequate and that we need to extend (\ref{log1}) and consider it as a matrix valued function. We remark that some progress on the topic of general L\'{e}vy-Khintchine formulae on symmetric spaces was obtained by different methods in the PhD thesis of Han Zhang \cite{Zhang}, however our results are more general and complete.

\section{Representation Theoretic Aspects of Generalised Eisenstein Integrals} \label{sec3}

\subsection{Generalised Eisenstein Integrals}

In the appendix, we describe the construction of the principal series representations of a connected semisimple Lie group $G$. The reader can find more background to this section there.
Define $M: = Z_{K}(\fa) = \{k \in K; \Ad(k)X = X~\mbox{for all}~X \in \fa\}$. Then $M$ is a compact Lie group whose Lie algebra is the centraliser of $\fa$ in $\fk$.
For a minimal parabolic subgroup $NAM$ of $G$, we choose  an irreducible unitary representation $\sigma$ of $M$ acting in a finite-dimensional vector space $\mathcal{H}_\sigma$ and $\lambda \in \fa^{*}_{\C}$.  The representation $1 \otimes \lambda \otimes \sigma$ of
 the minimal parabolic group $NAM$ may be induced to a representation $\xi_{\lambda}$ of $G = NAK$, which operates on $V_{K, \sigma}:=\{ f \in L^2(K, \mathcal{H}_\sigma): f(mk)= \sigma(m) f(k),~\mbox{for all}~m \in M~\mbox{and almost all}~k \in K\}$. This representation is directly related (see (\ref{kversion}) and preceding discussion for more details) to the following representation of $G$ on $L^2(K)$:

\begin{equation} \label{ind}
(\xi'_{\lambda}(g)f)(l) = e^{(\lambda + \rho)(A(lg))}f(u(lg)),
\end{equation}
for each $f \in V_{K, \sigma}, l \in K, g = n(g)\exp(A(g))u(g) \in G$.


We shall show in the appendix that for all $\lambda \in \at_{\C}, \xi_{\lambda}^{\prime}$ may be considered as a representation of $G$ in $L^2(K)$. If $\lambda$ is pure imaginary then the representation is unitary (see (\ref{kversion}) and the discussion that follows.) From now on we will take $\lambda$ to be pure imaginary and write it as $i\lambda$ where $\lambda \in \at$.  Notice that the restriction of $\xi'$ to $K$ is precisely the right regular representation of $K$ in $L^2(K)$. We will use this fact to decompose $\xi'$ into its $K$-types.

For each $\pi \in \widehat{K}$, let $V_{\pi}$ be the corresponding representation space so that for all $k \in K, \pi(k)$ acts unitarily in $V_{\pi}$. Note that $V_{\pi}$ is finite-dimensional for each $\pi \in \widehat{K}$. For each $u,v \in V_{\pi}, k \in K$ define $f_{u,v}^{\pi}(k):= \sqrt{d_{\pi}}\la \pi(k)u, v \ra$ and let ${\cal M}_{\pi}$ be the linear span of $\{f_{u,v}^{\pi}, u,v \in V_{\pi}\}$. Then ${\cal M}_{\pi}$ is a closed finite-dimensional subspace of $ L^{2}(K): = L^{2}(K, \mathbb{C})$ and by the Schur orthogonality relations, ${\cal M}_{\pi}$ and ${\cal M}_{\pi^{\prime}}$ are orthogonal if $\pi \neq \pi^{\prime}$. By the Peter-Weyl theorem,
$ L^{2}(K)= \bigoplus_{\pi \in \widehat{K}}{\cal M}_{\pi}$.

Recall that the mapping $v \mapsto v^*$ is a conjugate linear bijection between $V_{\pi}$ and $V_{\pi}^*$, where $v^*(u) = \langle u,v \rangle$ for $u,v \in V_{\pi}$, and that $V_{\pi}^*$ becomes an inner product space when we define $\langle v^*,u^* \rangle_{V_{\pi}^*}=\langle u,v \rangle_{V_{\pi}}$. Then if $\pi \in \hat{K}$, the {\it conjugate representation} is defined by $$\langle \bar{\pi}(k) u^*,v^* \rangle_{V_{\pi}^*}=\langle \pi(k^{-1}) v,u \rangle_{V_{\pi}} = \overline {\langle \pi (k) u,v \rangle_{V_\pi}}$$
for all $k \in K$.

For each $\pi \in \widehat{K}$, we define a linear mapping from $\gamma_{\pi}:V_{\pi} \otimes V_{\pi}^{*} \rightarrow {\cal M}_{\pi}$ by linear extension of the prescription
\begin{equation} \label{ten}
\gamma_{\pi}(u \otimes v^*) = f_{u,v}^{\pi},
\end{equation}
for each $u, v \in V_{\pi}$. The mapping $\gamma_{\pi}$ is clearly well-defined.

\begin{lemma} \label{ten1}
For each $\pi \in \widehat{K}, \gamma_{\pi}$ is a unitary isomorphism between $V_{\pi} \otimes V_{\pi}^{*}$ and ${\cal M}_{\pi}$. Moreover if $\pi_1, \pi_2 \in \widehat{K}$ with $\pi_1 \neq \pi_2$ then $\Ran(\gamma_{\pi_2}) \bot \Ran(\gamma_{\pi_1})$.
\end{lemma}
\noindent {\bf Proof.} It is clear from the definition (\ref{ten}) that $\gamma_{\pi}$ is surjective and hence bijective. The rest follows by Schur orthogonality. Indeed for $u_1, v_1 \in V_{\pi_1},  u_2, v_2 \in V_{\pi_2}$,
\bean \la \gamma_{\pi_1}(u_1 \otimes v_1^*), \gamma_{\pi_2}(u_2 \otimes v_2^*) \ra_{L^{2}(K)} & = & \int_{K}f_{u_1,v_1}^{\pi_1}(k)\overline{f_{u_2,v_2}^{\pi_2}}(k)dk \\
& = & \int_{K}\sqrt{d_{\pi_1}d_{\pi_2}} \la \pi_1(k)u_1, v_1 \ra_{V_{\pi_1}} \la v_2, \pi_2(k)u_2\ra_{V_{\pi_2}}dk \\
& = & \la u_1, u_2 \ra \la v_2,  v_1 \ra \delta_{\pi_1, \pi_2} \\
& = & \la u_1 \otimes v_1^*,  u_2 \otimes v_2^* \ra \delta_{\pi_1, \pi_2}. ~~~~~~~~~~~~~~\hfill \Box \eean

\vspace{5pt}


For each $\lambda \in \at, \pi_1, \pi_2 \in \widehat{K}, g \in G$ define $\Phi_{\lambda, \pi_1, \pi_2}(g) \in {\cal L}(V_{\pi_1} \otimes V_{\pi_2}^{*})$ by

\begin{equation} \label{supergen1}
\Phi_{\lambda,\pi_1, \pi_2}(g): = \sqrt{d_{\pi_{1}}d_{\pi_{2}}}\int_{K}e^{(i\lambda + \rho)(A(kg))}(\pi_1(u(kg)) \otimes \overline{\pi_2}(k))dk.
\end{equation}

We call $\Phi_{\lambda, \pi_1, \pi_2}$ a {\it generalised Eisenstein integral}. If we take $\pi_1$ to be the trivial representation acting in $\C$ then we may identify $\C \otimes V_{\pi_2}$ with $V_{\pi_2}$ and then (\ref{supergen1}) yields $ \Phi_{\lambda, \pi_2}(g):=\Phi_{\lambda, \pi_1, \pi_2}(g) = \sqrt{d_{\pi_{2}}}\int_{K}e^{(i\lambda + \rho)(A(kg))}\overline{\pi_2}(k)dk$ which is (a scalar multiple of) the Eisenstein integral (\ref{gsph}). If $\pi_2$ is also taken to be trivial then (\ref{supergen1}) reduces to the usual spherical function on non-compact semisimple Lie groups. Eisenstein integrals were used by van den Ban \cite{VB} in his analysis of the principal series representations on reductive groups; more recently, van den Ban and Schlichtkrull \cite{VBS} proved a full Plancherel decomposition for symmetric spaces using Eisenstein integrals. A lovely exposition of this work is given in \cite{VB2} (see also \cite{OS}).

\vspace{5pt}

Now we come to a key structural result which gives the connection between group representations and generalised Eisenstein integrals (note that here, and in the sequel, we always write $\xi$ instead of $\xi^{\prime}$ to keep the notation simple) :

\begin{theorem} \label{supergen}
For each $\lambda \in \at, \pi_1, \pi_2 \in \widehat{K}, g \in G, u_1, v_1 \in V_{\pi_1}, u_2, v_2 \in V_{\pi_2}$,
\begin{equation} \label{keystruct}
 \la \Phi_{\lambda,\pi_1, \pi_2}(g)(u_1 \otimes u_2^*), v_1 \otimes v_2^* \ra_{V_{\pi_1} \otimes V_{\pi_2}^{*}} = \la \xi_{\lambda}(g)\gamma_{\pi_1}(u_1 \otimes v_1^*), \gamma_{\pi_2}(u_2 \otimes v_2^*) \ra_{L^{2}(K)}
\end{equation}
\end{theorem}
\noindent {\bf Proof.} By (\ref{ten}), \bean & & \la \xi_{\lambda}(g)\gamma_{\pi_1}(u_1 \otimes v_1^*), \gamma_{\pi_2}(u_2 \otimes v_2^*) \ra_{L^{2}(K)}\\
& = & \la \xi_{\lambda}(g)f_{u_1,v_1}^{\pi_1}, f^{\pi_2}_{u_2,v_2} \ra_{L^{2}(K)} \\
& = & \int_{K}e^{(i\lambda + \rho)(A(kg))}f_{u_1,v_1}^{\pi_1}(u(kg))\overline{f^{\pi_2}_{u_2,v_2}(k)}dk \\
& = &  \sqrt{d_{\pi_1}d_{\pi_2}}\int_{K}e^{(i\lambda + \rho)(A(kg))} \la \pi_1(u(kg))u_1, v_1 \ra \overline{\la \pi_2(k)u_2, v_2 \ra}dk \\
& = & \sqrt{d_{\pi_1}d_{\pi_2}}\int_{K}e^{(i\lambda + \rho)(A(kg))} \la \pi_1(u(kg))u_1, v_1 \ra \la \overline{\pi_2}(k)u_2^*, v_2^* \ra dk \\
& = &  \sqrt{d_{\pi_1}d_{\pi_2}}\int_{K}e^{(i\lambda + \rho)(A(kg))} \la (\pi_1(u(kg)) \otimes \overline{\pi_2}(k))u_1 \otimes u_2^*, v_1
 \otimes v_2^* \ra dk \\
& = & \la \Phi_{\lambda, \pi_1, \pi_2}(g)(u_1 \otimes u_2^*), v_1 \otimes v_2^* \ra_{V_{\pi_1} \otimes V_{\pi_2}^{*}}~~~~~~~~~~~\hfill \Box\eean

\vspace{5pt}

We will want to work with matrix elements of $\Phi_{\lambda,\pi_1, \pi_2}(\cdot)$. To this effect let $(e_{\pi}^{i}, 1 \leq i \leq d_{\pi})$ be an orthonormal basis in $V_{\pi}$. For ease of notation we will always use $(h_{\pi}^{i}, 1 \leq i \leq d_{\pi})$ for the corresponding dual orthonormal basis in $V_{\pi}^{*}$ so that $h_{\pi}^{i}: = (e_{\pi}^{i})^{*}$ for all $1 \leq i \leq d_{\pi}$.  
Define for all $\lambda \in \at, \pi_1, \pi_2 \in \widehat{K}, g \in G, 1 \leq i,k \leq d_{\pi_1}, 1 \leq j,l \leq d_{\pi_2}$,
\begin{equation} \label{matent}
\Phi_{\lambda,\pi_1,\pi_2}^{(i,j), (k,l)}(g): = \la \Phi_{\lambda,\pi_1, \pi_2}(g)(e_{\pi_1}^{i} \otimes h_{\pi_2}^{j}), e_{\pi_1}^{k} \otimes h_{\pi_2}^{l} \ra,
\end{equation}
where we emphasise that the indices $i$ and $k$ are associated with the representation $\pi_1$ while $j$ and $l$ are associated with $\pi_2$. So we can identify the linear operator $\Phi_{\lambda,\pi_1, \pi_2}(g)$ acting on the finite dimensional vector space $V_{\pi_1} \otimes V_{\pi_2}^{*}$ with the $d_{\pi_1}d_{\pi_2} \times d_{\pi_1}d_{\pi_2}$ matrix $(\Phi_{\lambda,\pi_1,\pi_2}^{(i,j), (k,l)}(g))$ in the usual manner. It follows from (\ref{matent}) and (\ref{keystruct}) that the mapping $g \rightarrow \Phi_{\lambda,\pi_1, \pi_2}(g)$ from $G$ to $V_{\pi_1} \otimes V_{\pi_2}^{*}$ is $C^{\infty}$. We equip $V_{\pi_1} \otimes V_{\pi_2}^{*}$ with the usual Euclidean norm.

\begin{theorem} \label{contr} For all $\lambda \in \at, \pi_1, \pi_2 \in \widehat{K}, g \in G, \frac{1}{\sqrt{d_{\pi_{1}}d_{\pi_{2}}}}\Phi_{\lambda,\pi_1, \pi_2}(g)$ is a contraction in $V_{\pi_1} \otimes V_{\pi_2}^{*}$.
\end{theorem}
\noindent {\bf Proof.} Using (\ref{keystruct}), the Cauchy-Schwarz inequality and Lemma \ref{ten1}  we find that for all $1 \leq i,k \leq d_{\pi_1}, 1 \leq j,l \leq d_{\pi_2}, g \in G$,
\bean |\Phi_{\lambda,\pi_1,\pi_2}^{(i,j), (k,l)}(g)| & = & |\la \xi_{\lambda}(g)\gamma_{\pi_1}(e_{\pi_1}^{i} \otimes h_{\pi_2}^{j}), \gamma_{\pi_2}(e_{\pi_1}^{k} \otimes h_{\pi_2}^{l}) \ra|\\
& \leq & ||\gamma_{\pi_1}(e_{\pi_1}^{i} \otimes h_{\pi_2}^{j})||.||\gamma_{\pi_2}(e_{\pi_1}^{k} \otimes h_{\pi_2}^{l})||\\
& = & ||e_{\pi_1}^{i} \otimes h_{\pi_2}^{j}||.||e_{\pi_1}^{k} \otimes h_{\pi_2}^{l}|| = 1. \eean
So we have $\max_{1 \leq i,k \leq d_{\pi_1}, 1 \leq j,l \leq d_{\pi_2}}|\Phi_{\lambda,\pi_1,\pi_2}^{(i,j), (k,l)}(g)| \leq 1$.
The result then follows by a standard matrix estimate. $\hfill \Box$

\vspace{5pt}

The next result demonstrates that the group composition rule manifests as matrix multiplication of generalised Eisenstein integrals:

\begin{theorem} \label{matrixmult} For all $\lambda \in \at, \pi_1, \pi_2 \in \widehat{K}, g,h \in G, 1 \leq i,k \leq d_{\pi_1}, 1 \leq j,l \leq d_{\pi_2}$,
\begin{equation} \label{mmulteq}
\Phi_{\lambda,\pi_1, \pi_2}^{(i,j), (k,l)}(gh) = \sum_{\eta \in \widehat{K}}\sum_{p,q = 1}^{d_{\eta}}\Phi_{\lambda,\pi_1, \eta}^{(i,p), (k,q)}(h)\Phi_{\lambda,\eta, \pi_2}^{(p,j), (q,l)}(g) \end{equation}
\end{theorem}
\noindent {\bf Proof.} By (\ref{matent}) and Theorem \ref{supergen}
\bean \Phi_{\lambda,\pi_1, \pi_2}^{(i,j), (k,l)}(gh) & = &  \sqrt{d_{\pi_1}d_{\pi_2}}\la \xi(gh) (\la \pi_1(\cdot)e_{\pi_1}^{i},  e_{\pi_1}^{k}\ra), \la \pi_2(\cdot)e_{\pi_2}^{j}, e_{\pi_2}^{l} \ra \ra_{L^{2}(K)}\\
& = & \sqrt{d_{\pi_1}d_{\pi_2}} \la \xi(h)(\la \pi_1(\cdot)e_{\pi_1}^{i}, e_{\pi_1}^{k}\ra), \xi(g^{-1})(\la \pi_2(\cdot)e_{\pi_2}^{j}, e_{\pi_2}^{l} \ra)\ra_{L^{2}(K)} \\
& = &  \sqrt{d_{\pi_1}d_{\pi_2}}\sum_{\eta \in \widehat{K}}d_{\eta}\sum_{p,q =1}^{d_{\eta}}\la \xi(h)(\la \pi_1(\cdot)e_{\pi_1}^{i},e_{\pi_1}^{k}\ra), \la \eta(\cdot)e_{\eta}^{p}, e_{\eta}^{q} \ra \ra_{L^{2}(K)}\\ & \times &  \la \xi(g)(\la \eta(\cdot)e_{\eta}^{p}, e_{\eta}^{q} \ra), \la \pi_2(\cdot)e_{\pi_2}^{j}, e_{\pi_2}^{l} \ra \ra_{L^{2}(K)} \eean
and the result follows, again by Theorem \ref{supergen}. $\hfill \Box$

\vspace{5pt}

We can write (\ref{mmulteq}) succinctly as
\begin{equation} \label{mmulteq1}
\Phi_{\lambda,\pi_1, \pi_2}(gh) = \sum_{\eta \in \widehat{K}}\Phi_{\lambda,\pi_1, \eta}(h)\Phi_{\lambda,\eta, \pi_2}(g)
\end{equation}

Before we continue with our study of generalised Eisenstein integrals, we first recall that if $S$ is a countable set and $(H_{s}, s \in S)$ is a family of complex separable Hilbert spaces then we may form the direct sum $H: = \bigoplus_{s \in S}H_{s}$ comprising vectors $x = \bigoplus_{s \in S}x_{s}$ (where each $x_{s} \in H_{s}$) for which $\sum_{s \in S}||x_{s}||^{2} < \infty$. If $(T_{s}, s \in S)$ is a collection of linear operators with $T_{s} \in {\cal L}(H_{s})$ for all $s \in S$ and $\sup_{s \in S}||T_{s}|| < \infty$ then it is easy to verify that $T:= \bigoplus_{s \in S}T_{s} \in {\cal L}(H)$ where
$$ Tx:= \bigoplus_{s \in S}T_{s}x_{s}.$$
Note that if $T_{s}$ is a contraction for all $s \in S$ then so is $T$. More generally if each $T_{s}$ is a densely defined linear operator with domain $D_{s}$ then we may still define the direct sum operator $T$ as above on the dense domain given by the linear subspace $D$ of the vector space direct sum $\bigoplus_{s \in S}D_{s}$ comprising sequences of vectors that vanish for all but finitely many entries. In particular if $T$ and $V$ are two such operators and $V$ leaves $D$ invariant, then $TV$ is another such operator and
\begin{equation} \label{dsm}
TV = \bigoplus_{s \in S}T_{s}V_{s}.
\end{equation}

 Let ${\cal H}(K):=\bigoplus_{\pi_1, \pi_2 \in \widehat{K}}V_{\pi_1}\otimes V_{\pi_2}^{*}$ and define for all $\lambda \in \at, g \in G$ a linear operator $\Phi_{\lambda}(g)$ acting on the dense domain $V(K)$ of ${\cal H}(K)$ comprising sequences in $\bigoplus_{\pi_1, \pi_2\in \widehat{K}}V_{\pi_1} \otimes V_{\pi_2}^{*}$ which vanish for all but finitely many entries by
$$ \Phi_{\lambda}(g) = \bigoplus_{\pi_1, \pi_2 \in \widehat{K}}\Phi_{\lambda,\pi_1, \pi_2}(g).$$
Note that for all $g \in G, \Phi_{\lambda}(g)$ leaves $V(K)$ invariant.

Using (\ref{dsm}), We can now interpret (\ref{mmulteq1}) as the composition of two linear operators:
\begin{equation} \label{mm}
   \Phi_{\lambda}(gh) = \Phi_{\lambda}(h)\Phi_{\lambda}(g).
\end{equation}

\subsection{Infinitesimal Structure}

We will also need some information  about the differentials of the representations $\xi_{\lambda}$. To that end let $\theta$ be a Cartan involution on $\g$ and let $\g = \fk \oplus \fp$ be the corresponding Cartan
decomposition. Let $\Sigma_{+}$ be the set
of positive (with respect to a given lexicographic ordering) restricted
roots. Let $g_{\lambda}$ be the root space associated to the
restricted root $\lambda$ and define $X_{\lambda} \in g_{\lambda}$
by ad$(H)(X_{\lambda}) = \lambda(H)X_{\lambda}$ for all $H \in
\textbf{a}$. It follows from the Iwasawa
decomposition for semisimple Lie algebras (see e.g. \cite{Kn1}, Proposition 6.4.3, p.373) that we can decompose
each $X \in \fp$ as $X = H + Y_{-}$ where $H \in
\textbf{a}$ and $Y_{-} = \sum_{\lambda \in \Sigma_{+}}(X_{\lambda}
- \theta X_{-\lambda})$. Note that $Y_{-} \in \fn$ which is
the Lie algebra of $N$. We will also find it useful to introduce $Y_{+} =  \sum_{\lambda
\in \Sigma_{+}}(X_{\lambda} + \theta X_{-\lambda})$

\begin{lemma} \label{diffrep} Let $\lambda \in \at$. If $X \in \fk$, then
\begin{equation} \label{drift2}
d\xi_{\lambda}(X) = X
\end{equation}
and if $X \in \fp$, then
\begin{equation} \label{drift3}
(d\xi_{\lambda}(X)f)(l) = <i\lambda + \rho,\Pi_{{\mathfrak a}}(\mbox{Ad}(l)X)> f(l) + Y_{+}f(l),
\end{equation}
for $f\in C^\infty(K), l \in K$, where $\Pi_{\mathfrak{a}}$ is the orthogonal projection from
$\fp$ to $\mathfrak{a}$.
\end{lemma}
\noindent {\bf Proof.}~By (\ref{ind}) we have
\begin{equation} \label{indag}
(\xi_{\lambda}(\exp(tX)f)(l) = e^{(i\lambda +
\rho)(A(l\exp(tX)))}f(u(lk_{t}),
\end{equation}
for each $t \in \R, f \in L^{2}(K), l \in K$ and where we have used the
Iwasawa decomposition $\exp(tX) = n_{t}a_{t}u(k_{t})$.

We first suppose that $X \in \fk$ and so $\exp(tX) \in K$
for all $t \in \R$. Hence $A(l\exp(tX)) = 0$ for all $l \in K$.
Consequently (\ref{indag}) takes the form
$$ (\xi_{\lambda}(\exp(tX)f)(l) = f(l\exp(tX)),$$
and so (\ref{drift2}) follows when take $f \in C^{1}(K)$, differentiate with respect to $t$ and then put $t=0$.
Now suppose that $X \in \fp$. For each
$H \in \textbf{a}$ we define $Y_{++}^{(H)} = \sum_{\lambda \in
\Sigma_{+}}\lambda(H)(X_{\lambda} + \theta X_{-\lambda})$ and
$Y_{--}^{(H)} = \sum_{\lambda \in
\Sigma_{+}}\lambda(H)(X_{\lambda} - \theta X_{-\lambda})$. Note
that $Y_{+}, Y_{++}^{(H)} \in \fk$ and $Y_{--}^{(H)} \in
\fn$. Since each $\theta \g_{\lambda} = \g_{-\lambda}$ it
is easy to check that $[H, Y_{-}] = Y_{++}^{(H)}$ and $[H, Y_{+}]
= Y_{--}^{(H)}$ for each $H \in \textbf{a}$. It follows by the
Campbell-Baker-Hausdorff formula that for all $t \in \R$,
$$ \exp(tX) =
\exp{(tH)}\exp{(tY_{-})}\exp{\left(-\frac{1}{2}t^{2}Y_{++}^{(H)}\right)}\cdots,$$
from which we verify that
$\left.\ds\frac{d}{dt}A(\exp(tX))\right|_{t=0} = H$.

It follows that

\begin{equation} \label{firstp}
\left.\ds\frac{d}{dt}\xi_{\lambda}(\exp(tX))f(l)\right|_{t=0} =
<(i\lambda + \rho)\Pi_{\mathfrak{a}}(\mbox{Ad}(l)X)>f(l).
\end{equation}

We also have (by similar arguments) for $f \in C^{\infty}(K)$,

\begin{equation} \label{secondp}
\left.\frac{d}{dt}f(u(lk_{t}))\right|_{t=0} = Y_{+}f(l).
\end{equation}

(\ref{drift3}) then follows from  (\ref{firstp}) and
(\ref{secondp}). $\hfill \Box$

\vspace{5pt}

For each $\lambda \in \at, \pi_{1}, \pi_{2} \in \widehat{K}, X \in \g, 1 \leq i,j \leq d_{\pi_1}, 1 \leq k, l \leq d_{\pi_2}$, we define
\begin{equation} \label{diffmat}
  \rho_{\lambda, \pi_1, \pi_2}^{(i,j), (k,l)}(X): = \left.\frac{d}{dt}\Phi_{\lambda,\pi_1, \pi_2}^{(i,j), (k,l)}(\exp(tX))\right|_{t=0}
\end{equation}

Indeed it follows from (\ref{keystruct}) that $\rho_{\lambda, \pi_1, \pi_2}^{(i,j), (k,l)}(X)$ is well-defined and that
\begin{equation} \label{diffmat1}
 \rho_{\lambda, \pi_1, \pi_2}^{(i,j), (k,l)}(X) = \la d\xi_{\lambda}(X)\gamma_{\pi_1}(e_{i}^{\pi_1} \otimes h_{j}^{\pi_1}), \gamma_{\pi_2}(e_{k}^{\pi_2} \otimes h_{l}^{\pi_2}) \ra.
\end{equation}

We will require the matrix $\rho_{\lambda, \pi_1, \pi_2}: = (\rho_{\lambda, \pi_1, \pi_2}^{(i,j), (k,l)}(X))$ acting on $V_{\pi_1} \otimes V_{\pi_2}^{*}$ and also the direct sum operator
$\rho_{\lambda}(X) : = \bigoplus_{\pi_1, \pi_2 \in \widehat{K}}\rho_{\lambda, \pi_1, \pi_2}(X)$ acting on the dense domain $V(K)$. Note that operators of the form $\rho_{\lambda}(X)\Phi_{\lambda}(g)$ and $\rho_{\lambda}(X)\rho_{\lambda}(Y)\Phi_{\lambda}(g)$  are well-defined on the domain $V(K)$ for all $\lambda \in \at, g \in G, X, Y\in \g$ and the results of the following lemma implicitly utilise this action.

\begin{lemma} \label{diffGan}
For all $\lambda \in \at, \pi_{1}, \pi_{2} \in \widehat{K},g \in G, X, Y \in \g$
\begin{equation} \label{diffGan1}
X\Phi_{\lambda, \pi_{1}, \pi_{2}}(g) = \sum_{\pi \in \widehat{K}}\rho_{\lambda, \pi_{1}, \pi}(X)\Phi_{\lambda, \pi, \pi_{2}}(g),
\end{equation}
\begin{equation} \label{diffGan2}
XY\Phi_{\lambda, \pi_{1}, \pi_{2}}(g) = \sum_{\pi \in \widehat{K}}\sum_{\eta \in \widehat{K}}\rho_{\lambda, \pi_{1}, \pi}(Y)\rho_{\lambda, \pi, \eta}(X)\Phi_{\lambda, \eta, \pi_{2}}(g).
\end{equation}

\end{lemma}
\noindent {\bf Proof.} For all $1 \leq i,j \leq d_{\pi_{1}}, 1 \leq k, l \leq d_{\pi_{2}}$, using (\ref{matent}) and (\ref{keystruct}) we have
\bean X\Phi_{\lambda, \pi_{1}, \pi_{2}}^{(i,j), (k,l)}(g) & = & \left.\frac{d}{dt}\Phi_{\lambda, \pi_{1}, \pi_{2}}^{(i,j), (k,l)}(g\exp(tX))\right|_{t=0}\\
& = & \left.\frac{d}{dt}\la \xi_{\lambda}(\exp(tX))f^{\pi_{1}}_{e_{i}^{\pi_{1}}, e_{j}^{\pi_{1}}},  \xi_{\lambda}(g^{-1})f^{\pi_{2}}_{e_{k}^{\pi_{2}}, e_{l}^{\pi_{2}}} \ra \right|_{t=0}\\
& = & \la d\xi_{\lambda}(X)f^{\pi_{1}}_{e_{i}^{\pi_{1}}, e_{j}^{\pi_{1}}},  \xi_{\lambda}(g^{-1})f^{\pi_{2}}_{e_{k}^{\pi_{2}}, e_{l}^{\pi_{2}}} \ra \eean
and the first result follows by similar computations to those that feature in the proof of Theorem \ref{matrixmult}. The second result is established by iterating this argument. $\hfill \Box$

\vspace{5pt}
The precise expressions for the coefficients $\rho_{\lambda, \pi_{1}, \pi}(X)$, for $X \in \mathfrak{g} = \mathfrak{k} + \mathfrak{p}$ are given in Lemma \ref{diffrep}.

\vspace{5pt}
It is clear from the argument in the proof that for all $X,Y, Z \in \g$ the mappings $g \rightarrow X\Phi_{\lambda, \pi_{1}, \pi_{2}}(g)$ and $g \rightarrow
YZ\Phi_{\lambda, \pi_{1}, \pi_{2}}(g)$ are bounded and continuous.

\section{Infinite Matrices and the Eisenstein Transform}\label{sec4}

\subsection{The Exponential of an Infinite Matrix}

In the next section we will need the exponential of an infinite matrix that is indexed by $\widehat{K} \times \widehat{K}$. First we introduce an identity matrix
$\delta_{\pi_{1}, \pi_{2}} = \left\{\begin{array}{c c} I &~\mbox{if}~\pi_{1} = \pi_{2}\\ 0 &~\mbox{if}~\pi_{1} \neq \pi_{2}\end{array} \right.$.

Let $\Upsilon(t) = (\Upsilon_{\pi_{1}, \pi_{2}}(t), \pi_{1}, \pi_{2} \in \widehat{K})$ be a family of infinite matrices (acting on $V(K)$) indexed by $t \geq 0$ which have the following properties:
\begin{enumerate}
\item[E(i)] $\Upsilon(s+t) = \Upsilon(s)\Upsilon(t)$ for all $s,t \geq 0$,
\item[E(ii)] $\Upsilon_{\pi_{1}, \pi_{2}}(0) = \delta_{\pi_{1}, \pi_{2}}$,
\item[E(iii)] The mapping $t \rightarrow \Upsilon_{\pi_{1}, \pi_{2}}(t)$ is continuous for all $\pi_{1}, \pi_{2} \in \widehat{K}$.
\item[E(iv)] The mapping $t \rightarrow \Upsilon_{\pi_{1}, \pi_{2}}(t)$ is differentiable at $t=0$ and there exists an infinite matrix $\Theta$ so that $\left.\frac{d}{dt}\Upsilon_{\pi_{1}, \pi_{2}}(t)\right|_{t=0} = \Theta_{\pi_{1}, \pi_{2}}$ for all $\pi_{1}, \pi_{2} \in \widehat{K}$.
\end{enumerate}

In this case we write $\Upsilon(t): = \Exp(t \Theta)$ and call $\Exp$ an ``infinite matrix exponential''. It clearly coincides with the usual exponential for finite dimensional matrices.  We also write $A: = \Exp(\Theta)$ whenever $A = \Upsilon(1)$ for some $(\Upsilon(t), t \geq 0)$ satisfying (E(i)) to (E(iv)).

Note that for E(iii) and E(iv) (and also in calculations to follow) we assume that some fixed norm is chosen on
each space of $d_{\pi_{1}}d_{\pi_{2}} \times d_{\pi_{1}}d_{\pi_{2}}$ matrices.

\subsection{The Eisenstein Transform}

We generalise the spherical transform of \cite{Gang1} and for each $\lambda \in \at$ we introduce the {\it Eisenstein transform} $\widehat{\mu}_{\lambda}$ of a Borel probability measure $\mu$ to be the matrix valued integral defined for each $\pi,\pi_2 \in \widehat{K}$ by:

\begin{equation} \label{Eis}
\widehat{\mu}_{\lambda, \pi_1, \pi_2} : = \int_{G}\Phi_{\lambda,  \pi_1, \pi_2}^{\#}(g)\mu(dg)
\end{equation}
for all $\lambda \in \at$.
This is related in an obvious way to the group Fourier transform defined in section \ref{sec74}. Now let $\mu_{1}$ and $\mu_{2}$ be Borel probability measures on $G$ and define their convolution $\mu_{1} * \mu_{2}$ in the usual way:
\begin{equation} \label{conv}
\int_{G}f(g)(\mu_{1} * \mu_{2})(dg) = \int_{G}\int_{G}f(gh) \mu_{1}(dg)\mu_{2}(dh),
\end{equation}
where $f$ is an arbitrary bounded Borel measurable function on $G$. Before we investigate the Eisenstein transform of the convolution of two measures, we state a useful technical result.

\begin{lemma} \label{Planch}
Let $U$ be a unitary representation of a Lie group $\Gamma$ acting on a complex separable Hilbert space $H$ and let $(e_{n}, \nN)$ be a complete orthonormal basis for $H$. Then for all $x, y \in H$ the infinite series $\sum_{n=1}^{\infty}\la U(g)x, e_{n} \ra \la e_{n}, U(h)y \ra$ converges uniformly in $g$ and $h$ to $\la U(g)x, U(h)y \ra$.
\end{lemma}
\noindent {\bf Proof.} The result follows from unitarity and the proof of the usual Parseval identity in Hilbert space. $\hfill \Box$

\begin{theorem} \label{convhatthm}
For each $\lambda \in \at$ (in the sense of matrix multiplication)
\begin{equation} \label{convhat}
(\widehat{\mu^{(1)} * \mu^{(2)}})_{\lambda} =   \widehat{\mu_{\lambda}^{(1)}}\widehat{\mu_{\lambda}^{(2)}},
\end{equation}
\end{theorem}

{\bf Proof.} From (\ref{mm}) and (\ref{conv}) we deduce that for each $\pi_{1}, \pi_{2} \in \widehat{K}, \lambda \in \at$,
\bean (\widehat{\mu^{(1)} * \mu^{(2)}})_{\lambda} & = & \int_{G}\int_{G}\Phi_{\lambda,  \pi_1, \pi_2}^{\#}(gh)\mu_{1}(dg)\mu_{2}(dh)\\
& = & \int_{G}\int_{G}\left(\sum_{\pi \in \widehat{K}}\Phi_{\lambda,\pi_1, \pi}^{\#}(g)\Phi_{\lambda,\pi, \pi_{2}}^{\#}(h)\right)\mu_{1}(dg)\mu_{2}(dh). \eean
The result follows on interchanging the infinite sum and the integrals. This is justified as follows. Let $S \subseteq \widehat{K}$ and define for all $\lambda \in \at, g,h \in G$
$$ M^{S}_{\lambda, \pi_{1}, \pi_{2}}(g,h) := \sum_{\pi \in S}\Phi^{\#}_{\lambda, \pi_{1}, \pi}(g)\Phi^{\#}_{\lambda, \pi, \pi_{2}}(h).$$
Then by Theorem \ref{matrixmult}
\bean M^{S}_{\lambda, \pi_{1}, \pi_{2}}(g,h) & = & \sum_{\pi \in \widehat{K}}\Phi^{\#}_{\lambda, \pi_{1}, \pi}(g)\Phi^{\#}_{\lambda, \pi, \pi_{2}}(h) - \sum_{\pi \in S^{c}}\Phi^{\#}_{\lambda, \pi_{1}, \pi}(g)\Phi^{\#}_{\lambda, \pi, \pi_{2}}(h)\\
& = & \Phi^{\#}_{\lambda, \pi_{1}, \pi_{2}}(gh) - M^{S^{c}}_{\lambda, \pi_{1}, \pi_{2}}(g,h). \eean

By the construction in the proof of Theorem \ref{matrixmult} and Lemma \ref{Planch} we can assert that given any $\epsilon > 0$ there exists a finite subset $S_{0}$ of $\widehat{K}$ so that if $S\subseteq \widehat{K}$ is any other finite subset with $S_{0} \subset S$ then $\sup_{g, h \in G}|| M^{S^{c}}_{\lambda,\pi_{1}, \pi_{2}}(\sigma,\tau)|| < \epsilon$.
But then by Theorem \ref{contr}
\bean \sup_{g, h \in G}|| M^{S}_{\lambda, \pi_{1}, \pi_{2}}(g,h)|| & \leq &   \sup_{g,h \in G}||\Phi^{\#}_{\lambda, \pi_{1}, \pi_{2}}(gh)|| + \sup_{g, h \in G}|| M^{S^{c}}_{\lambda, \pi_{1}, \pi_{2}}(g,h)||\\
& \leq & C_{\pi_{1}, \pi_{2}}\sqrt{d_{\pi_{1}}d_{\pi_{2}}} + \epsilon, \eean
where $C_{\pi_{1}, \pi_{2}} \geq 0$ depends on the choice of matrix norm.
The right hand side of the inequality is constant and thus $\nu$-integrable and so dominated convergence can be applied to interchange the sum and the integral.
$\hfill \Box$

\vspace{5pt}

{\bf Remark.} It is a difficult problem to find conditions for the injectivity of the mapping $\mu \rightarrow \widehat{\mu}$ beyond the known result when $\mu$ is $K$-bi-invariant (see Theorem 4.1 in \cite{Gang1}).  Indeed this will not be true even in the case where $\mu$ has an $L^{2}$-density . Some progress on the problem when $\mu$ has a $C^\infty$ density can be found in the generalised Paley Wiener theorem that is presented as Theorem 4.1 in \cite{Art}.

For an injective Fourier transform on $L^2(G)$, one needs to analyse the discrete series and all generalized principal series, that is representations of the form
$(\xi \otimes \lambda \otimes 1)\uparrow_Q^G,$ with $Q$ a cuspidal parabolic subgroup of $G,$ $\xi$ a discrete
series rep of $M_Q,$ and $\lambda$ a unitary character of $A_Q.$  In this paper, we have only dealt with the case where $Q$ is minimal. We defer investigation of the other series of representations to future work. We would like to thank Erik van den Ban for useful discussions on this point \cite{VB3}.

\section{The Extended Gangolli-L\'{e}vy Khintchine Formula}\label{sec5}

In this section we will need to make sense of expressions of the form ${\cal L}\Phi_{\lambda, \pi_{1}, \pi_{2}}$. Note that the matrix elements of $\Phi_{\lambda, \pi_{1}, \pi_{2}}$ do not necessarily vanish at infinity. We use the facts that the generalised Eisenstein integrals are bounded and have second order bounded derivatives (in the sense of Lie algebra actions) and that $(T_{t}, t \geq 0)$ extends to a locally uniformly continuous semigroup on the space $C_{b}(G)$ of bounded continuous functions on $G$ (see Lemma 3.1 in \cite{Sch}). Now define the space $C_{b}^{2}(G)$ in exactly the same way as $C^{2}(G)$ but with $C_{0}(G)$ replaced by $C_{b}(G)$ wherever the first space appears within the definition. Then ${\cal L}$ extends to a linear operator in $C_{b}(G)$ with $C_{b}^{2}(G) \subseteq \Dom({\cal L})$ and ${\cal L}\Phi_{\lambda, \pi_{1}, \pi_{2}}$ has meaning as the $d_{\pi_{1}}d_{\pi_{2}} \times d_{\pi_{1}}d_{\pi_{2}}$ matrix with entries ${\cal L}\Phi_{\lambda, \pi_{1}, \pi_{2}}^{(i,j),(k,l)}$.

We return to the study of L\'{e}vy processes. Instead of
(\ref{log1}), we will find it more profitable to consider

\begin{equation} \label{log2}
\alpha_{\lambda, \pi_{1}, \pi_{2}}(t) :=\E(\Phi^{\#}_{\lambda, \pi_{1},
\pi_{2}}(Z(t))) = \int_{G}\Phi_{\lambda, \pi_{1},
\pi_{2}}(g^{-1})\mu_{t}(dg),
\end{equation}
 for all $\lambda \in \at, \pi_1, \pi_2 \in \widehat{K}, t \geq 0$ as an infinite matrix. Let us write this matrix as $\alpha_{\lambda}(t)$. Our goal will be to find another infinite matrix $\psi_{\lambda}$ with generic matric entry $\psi_{\lambda, \pi_{1}, \pi_{2}}$ which has the property that for all $t \geq 0$:

\begin{equation} \label{log3}
\alpha_{\lambda }(t) = \Exp(t \psi_{\lambda}).
\end{equation}


\begin{theorem} \label{EGLK2} Let $(\mu_{t}, t \geq 0)$ be the convolution semigroup of probability measures generated by ${\cal L}$. Then (\ref{log3}) holds. Furthermore for each $\lambda \in  \fa^{*}_{\C}$

\begin{equation} \label{symb}
\psi_{\lambda} = -b^{i}\rho_{\lambda}(X_{i}) + a^{ij}\rho_{\lambda}(X_{i})\rho_{\lambda}(X_{j}) + \eta_{\lambda},
\end{equation}
where $\eta_{\lambda}$ is the matrix whose diagonal entries are
$$\int_{G}(\Phi^{\#}_{\lambda, \pi, \pi}(\tau) - 1 +
x^{i}(\tau)\rho_{\lambda,\pi, \pi}(X_{i})) \nu(d\tau)$$ and
off-diagonal entries are $$\int_{G}(\Phi^{\#}_{\lambda, \pi_1,
\pi_2}(\tau) + x^{i}(\tau)\rho_{\lambda,\pi_1,\pi_2}(X_{i})) \nu(d\tau).$$
\end{theorem}
\noindent {\bf Proof.} We first check that the integrals defining the entries of the infinite matrix $\eta_{\lambda}$ in (\ref{symb})
are finite. To establish this we must show that the integrand is $O\left(\sum_{i=1}^{n}x_{i}^{2}\right)$ in a canonical co-ordinate neighbourhood $V$ of $e$.
Note that we may choose $x_{1}, \ldots, x_{n}$ so that $x(g^{-1}) = -x(g)$ for all $g \in V$. For each $g \in V$, we write $g = \exp(X)$ where $X \in \g$ and we use the fact that for each $\lambda \in \at, \xi_{\lambda}(\exp(X)) = \exp(d\pi(\xi_{\lambda}(X)))$ where the second $\exp$ is understood in the sense of functional calculus for skew-adjoint operators. Now using a Taylor expansion, (\ref{keystruct}), (\ref{diffmat1}) and Lemma \ref{diffGan} we find that for all $\pi_{1}, \pi_{2} \in \widehat{K}, g \in V$:
\bean \Phi_{\lambda, \pi_{1}, \pi_{2}}(g^{-1}) & = & \delta_{\pi_{1}, \pi_{2}} - \sum_{i=1}^{n}x_{i}(g)\rho_{\lambda, \pi_{1}, \pi_{2}}(X_{i})\\
& + & \sum_{i,j=1}^{n}x_{i}(g)x_{j}(g)\sum_{\pi, \eta \in \widehat{K}}\rho_{\lambda, \pi_{1}, \pi}(X_{i})\rho_{\lambda, \pi, \eta}(X_{j})\Phi_{\lambda, \eta, \pi_{2}}(g^{\prime}), \eean where $g^{\prime} \in V$, from which the required integrability follows.

To prove the theorem we need to show that E(i) to E(iv) are satisfied for $(\alpha_{\lambda}(t), t \geq 0)$. E(i) follows from Theorem \ref{convhatthm} since we have $\mu_{t+s} = \mu_{t} * \mu_{s}$ for all $s, t \geq 0$, E(ii) holds since $\mu_{0} = \delta_{e}$ and $\Phi_{\lambda, \pi_{1}, \pi_{2}}(e) = \delta_{\pi_{1}, \pi_{2}}$ for all $\pi_{1}, \pi_{2} \in \widehat{K}$. E(iii) is clear from the fact that $t \rightarrow \mu_{t}$ is weakly continuous. To establish E(iv) we argue as in the Remark on page 23 of \cite{Liao} (see also equation (8) in \cite{LW}) to see that
$$ \left.\frac{d}{dt}\int_{G}\Phi_{\lambda}^{\#}(g)\mu_{t}(dg)\right|_{t=0} = {\cal L}\Phi_{\lambda}^{\#}(e) = \psi_{\lambda},$$
where the last identity is obtained by explicit calculation. $\hfill \Box$

\begin{cor} Let $\mu$ be an infinitely divisible probability measure on $G$. Then for all $\lambda \in \at$ we have
$$ \int_{G}\Phi^{\#}_{\lambda}(g)\mu(dg) = \Exp(\psi_{\lambda}),$$
where the matrix $\psi_{\lambda}$ is as in Theorem \ref{EGLK2}.
\end{cor}
\noindent {\bf Proof.}~From the remarks at the beginning of section 6 in \cite{DM} we see that any connected semisimple Lie group has the ``embedding property'' that $\mu = \mu_{1}$ for some convolution semigroup of probability measures $(\mu_{t}, t \geq 0)$ on $G$. The result then follows from Theorem \ref{EGLK2}. $\hfill \Box$

\vspace{5pt}

We have a partial converse to this result, but the last part should be read in conjunction with the Remark at the end of section 4.

\begin{prop} Let $\psi_{\lambda}^{(b,a,\nu)}$ be an infinite matrix of the form given in Theorem \ref{EGLK2} (where we emphasise the dependence on the triple $(b, a, \nu)$ comprising a vector, a non-negative definite matrix $a$ and a L\'{e}vy measure $\nu$ on $G-\{e\}$.) Then there is an infinitely divisible probability measure $\mu$ on $G$ so that \begin{equation} \label{Eisnew1} \int_{G}\Phi^{\#}_{\lambda}(g)\mu(dg) = \Exp(\psi_{\lambda}) \end{equation} for all $\lambda \in \at$. Furthermore if $\mu$ is such that the Eisenstein transform $\mu \rightarrow \widehat{\mu}$ is injective then $\mu$ is uniquely defined by (\ref{Eisnew1}). \end{prop}

{\it Proof.}  Using the triple $(b, a, \nu)$ we construct a linear operator of the form (2.6) acting on $C^{2}(G)$.  By Theorem \ref{Hunt} this operator extends to be the infinitesimal generator of a unique convolution semigroup of probability measures $(\mu_{t} , t \geq 0)$. Let $(T_{t}, t \geq 0)$ be the corresponding semigroup of convolution operators. Then by the procedure of Theorems 4.1 and 4.2 we have for all $t \geq 0$.
$$  T_{t}\Phi_{\lambda, \pi_{1}, \pi_{2}}(e) = \int_{G}\Phi_{\lambda, \pi_{1}, \pi_{2}}(g^{-1})\mu_{t}(dg) = \Exp(t \psi_{\lambda,\pi_{1}, \pi_{2}}^{(b,a,\nu)})$$ and $\mu:=\mu_{1}$.

Now suppose that we can find another triple $(b^{\prime}, a^{\prime}, \nu^{\prime})$ so that $\psi_{\lambda}^{(b,a,\nu)} = \psi_{\lambda}^{(b^{\prime}, a^{\prime}, \nu^{\prime})}$. Then by the above construction we can find a convolution semigroup $(\mu_{t}^{\prime} , t \geq 0)$ for which
$\int_{G}\Phi_{\lambda, \pi_{1}, \pi_{2}}(g^{-1})\mu_{t}(dg) = \int_{G}\Phi_{\lambda, \pi_{1}, \pi_{2}}(g^{-1})\mu_{t^{\prime}}(dg)$ and so if the Eisenstein transform is injective then $\mu_{t} = \mu_{t^{\prime}}$ for all $t \geq 0$.  $\hfill \Box$

\section{Extension to Symmetric Spaces}\label{sec6}

Let $\mu$ be a right $K$-invariant Borel probability measure on $G$ so that $\mu(Ak) = \mu(A)$ for all $k \in K$ and all $A \in {\cal B}(G)$.  We now consider the Eisenstein transform $\widehat{\mu}_{\lambda}$ of $\mu$ as defined in (\ref{Eis}). Let $\delta$ denote the trivial representation of $K$.

\begin{theorem} \label{trivgood} If $\mu$ is a right $K$-invariant probability measure on $G$ then for all $\lambda \in \at, \pi \in \widehat{K},\widehat{\mu}_{\lambda, \pi_1, \pi} = 0$ except when $\pi_1 = \delta$.
\end{theorem}
\noindent{\bf Proof.} By (\ref{Eis}), (\ref{keystruct}) and using the right $K$-invariance of $\mu$ we obtain for all $\pi_1 \in \widehat{K}, u_1, v_1 \in V_{\pi_1}, u, v , \in V_{\pi}$,
\bean \la \widehat{\mu}_{\lambda, \pi_1, \pi}(u_1 \otimes u^{*}), v_1 \otimes v^{*} \ra & = & \int_{G} \la \Phi^{\#}_{\lambda, \pi_1, \pi}(g)(u_1 \otimes u^{*}), v_1 \otimes v^{*} \ra \mu(dg)\\
& = & \int_{G} \la \xi_{\lambda}(g)f^{\pi_1}_{u_1, v_1}, f^{\pi}_{u,v} \ra \mu(dg)\\
& = &  \int_{K}\int_{G} \la \xi_{\lambda}(g)f^{\pi_1}_{u_1, v_1}, f^{\pi}_{u,v} \ra \mu(dg)dk\\
& = & \int_{K}\int_{G} \la \xi_{\lambda}(gk)f^{\pi_1}_{u_1, v_1}, f^{\pi}_{u,v} \ra \mu(dg)dk\\
& = & \int_{G} \left \la \left(\int_{K}\xi_{\lambda}(k)f^{\pi_1}_{u_1, v_1}dk\right), \xi(g^{-1})f^{\pi}_{u,v} \right\ra \mu(dg)\eean
by Fubini's theorem.
However by Peter-Weyl theory, for each $l \in K$
$$ \int_{K}\xi_{\lambda}(k)f^{\pi_1}_{u_1, v_1}(l)dk = \int_{K}\la \pi_1(l)\pi_1(k)u_1, v_1\ra dk = \int_{K} \la \pi_1(k)u_1, v_1\ra dk = 0 $$
unless $\pi_1 = \delta$, and the result follows. $\hfill \Box$

\vspace{5pt}

Now let $\mu$ be right $K$-invariant. Then by Theorem \ref{trivgood}, the non-trivial values of $\widehat{\mu}_{\lambda}$ are given by
$$\widehat{\mu}_{\lambda, \pi}:= \widehat{\mu}_{\lambda, \delta, \pi} = \int_{G}\Phi^{\#}_{\lambda, \delta, \pi}(g)\mu(dg) = \int_{G}\Phi^{\#}_{\lambda, \pi}(g)\mu(dg),$$
where we recall that $\Phi_{\lambda, \pi}$ is the Eisenstein integral (\ref{gsph}). For the statement of the next theorem we will also require the (Harish-Chandra) spherical function, $\Phi_{\lambda}:= \Phi_{\lambda, \delta}$.

Now let $M:=G/K$ be the Riemannian globally symmetric space whose points are the left cosets $\{gK, g \in G\}$ and let $\natural: G \rightarrow G/K$ be the canonical surjection. Let $C_{K}(G)$ denote the closed subspace of $C_{0}(G)$ comprising those functions that are right $K$-invariant and observe that there is a natural isometric isomorphism between $C_{0}(G/K)$ and $C_{K}(G)$ given by the prescription $f \rightarrow f \circ \natural$ for each $f \in C_{0}(G/K)$.  Let $\mu$ be a right $K$-invariant Borel probability measure on $G$. Then $\widetilde{\mu} = \mu \circ \natural^{-1}$ defines a Borel probability measure on $M$. We define the {\it Eisenstein transform} of $\tilde{\mu}$ to be $\widehat{\mu}$ and note that this is well-defined as the Eisenstein transform of $\mu$ is invariant under the right $K$ action on the measure.

Now let $(\mu_{t}, t \geq 0)$ be a weakly continuous convolution semigroup of right invariant probability measures on $G$ where instead of requiring $\mu_{0} = \delta_{e}$ (which is not right invariant) we take $\mu_{0}$ to be normalised Haar measure on $K$ (equivalently $\mu_{0}(\cdot) = \int_{K}\delta_{e}(\cdot k)dk$, in the weak sense.) This may arise as the law of a $G$-valued L\'{e}vy process $Z$ with right $K$-invariant laws for which $Z(0)$ is uniformly distributed on $K$ and is independent of $Z(t)$ for $t > 0$. The measures $(\mu_{t}, t \geq 0)$ induce a $C_{0}$ contraction semigroup $(T_{t}, t \geq 0)$ having generator ${\cal L}$ on $C_{0}(G/K)$ in the usual manner, i.e. $T_{t}f = f * \mu_{t}$ for all $f \in C_{0}(G/K), t \geq 0$. We induce a $C_{0}$ semigroup of operators $(\widetilde{T_{t}}, t \geq 0)$ on $C_{0}(G/K)$ by the prescription
$$ (\widetilde{T_{t}}f) \circ \natural = T_{t}(f \circ \natural),$$ for all $t \geq 0, f \in C_{0}(G/K)$. Let $\widetilde{L}$ be the infinitesimal generator of $(\widetilde{T_{t}}, t \geq 0)$ then $(\widetilde{L}f) \circ \natural = L(f \circ \natural)$ for all $f \in \Dom(\widetilde{L}): = \{f \in C_{0}(G/K); f \circ \natural \in \Dom(L)\}$. We define $C_{2}(G/K): = \{f \in C_{0}(G/K); f \circ \natural \in C_{2}(G)\}$. Then $C_{2}(G/K)$ is dense in $C_{0}(G/K)$ and $C_{2}(G/K) \subseteq \Dom(\widetilde{L})$.
By restricting the matrix calculation in the proof of Theorem \ref{EGLK2} to preserve the row containing Eisenstein integrals we deduce the following which can be seen as a generalised L\'{e}vy-Khintchine formula on the symmetric space $M$. In the following, for each $\lambda \in \at, t \geq 0$ we regard $\widehat{\mu_{t}}(\lambda)$ and $\psi_{\lambda}$ as row-vectors where the value of the former at $\pi \in \widehat{K}$ is $\widehat{\mu_{t}}(\lambda, \pi): = \widehat{\mu_{t}}(\lambda,\delta, \pi)$.

\begin{theorem} \label{EGLK3} Let $(\mu_{t}, t \geq 0)$ be a right $K$-invariant convolution semigroup of probability measures on $G$. Then for all $\lambda \in \at, t > 0$ we have that $\widehat{\mu_{t}}(\lambda)$ is the row vector of $\Exp(t\psi_{\lambda})$ whose entries are of the form $(\delta, \pi)$ for $\pi \in \widehat{K}$. Furthermore for all $\pi \in \widehat{K}$
\begin{equation} \label{LKMan}
       \left.\frac{d}{dt}\widehat{\mu_{t}}(\lambda, \pi)\right|_{t=0} = \psi_{\lambda, \pi}
\end{equation}
where
\bean \psi_{\lambda, \pi} &  = & \psi_{\lambda, \delta, \pi}\\
& = & -b^{i}\rho_{\lambda, \delta, \pi}(X_{i}) + \sum_{\eta \in \widehat{K}}a^{ij}\rho_{\lambda, \delta, \eta}(X_{i})\rho_{\lambda, \eta, \pi}(X_{j}) + \eta_{\lambda, \pi}. \eean

Here $\eta_{\lambda, \pi}$ is the row vector whose $(\delta, \delta)$ entry is
$$\int_{G}(\Phi^{\#}_{\lambda}(\tau) - 1) \nu(d\tau)$$ and
$(\delta, \pi)$ entries (for $\pi \neq \delta$) are $$\int_{G}(\Phi^{\#}_{\lambda,\pi}(\tau) + x^{i}(\tau)\rho_{\lambda,\delta, \pi}(X_{i}))\nu(d\tau).$$
\end{theorem}

\noindent{\bf Proof.} This follows from the proof of Theorem \ref{EGLK2}. $\hfill \Box$

\vspace{5pt}

In particular, when $\mu_{t}$ is $K$-bi-invariant for all $t \geq 0$ then we recover the characteristic exponent in Gangolli's L\'{e}vy-Khintchine formula (for which see \cite{Gang1}, \cite{LW}, \cite{App0}) as $\psi_{\lambda, \delta}$.


\section{Example: $SU(1,1)$ and Hyperbolic Space}\label{sec7}
The group $SU(1,1)$ is by definition the set of $2 \times 2$ complex matrices of the form $\left(
                                                                                 \begin{array}{cc}
                                                                                   a & b \\
                                                                                   \bar{b} & \bar{a}\\
                                                                                 \end{array}
                                                                               \right) $ for which $|a|^2 -|b|^2 = 1.$
The group acts on the unit disc $\mathcal{D}$ in $\C^2$ by fractional linear transformations.

 It is worth noting that the Cayley transform $z\mapsto \ds\frac{z-i}{z+i}$ takes the upper half plane to $\mathcal{D}$, and that conjugation with the corresponding matrix
$\left(
   \begin{array}{cc}
     1 & -i \\
     1 & i \\
   \end{array}
 \right)$
gives an isomorphism between $SU(1,1)$ and $SL(2, \R)$. Thus the representation theory of $SU(1,1)$, which we are about to describe, is the representation theory of $SL(2, \R)$ in another guise.

The Lie algebra $\mathfrak{su}(1,1)$ is given by the $2 \times 2$ complex matrices of the form
$\left(
   \begin{array}{cc}
     ix & b \\
     \bar{b} & -ix \\
   \end{array}
 \right)$
 where $x \in \R$ and $b \in \C$.

 The Cartan involution is  $\theta(A) = -\bar{A}^t$ for $A \in \mathfrak{su}(1,1)$. The fixed points of $\theta$ are the real linear span $\mathfrak{k}$ of the matrix
$ \left(
   \begin{array}{cc}
     i & 0 \\
     0 & -i \\
   \end{array}
 \right)$
 which corresponds to the subgroup $SO(2)=\left\{k(\psi)=\left(
                                           \begin{array}{cc}
                                             \cos \psi &  \sin \psi \\
                                             -\sin \psi & \cos \psi \\
                                           \end{array}
                                         \right): \psi \in [0, 2 \pi)\right\}$.

 Letting $\mathfrak{p}$ be the linear span of $\left\{\left(
                                 \begin{array}{cc}
                                   0 & 1 \\
                                   1 & 0 \\
                                 \end{array}
                               \right),
\left(
  \begin{array}{cc}
    0 & i \\
    -i & 0 \\
  \end{array}
\right)\right \}$, we have $\mathfrak{su}(1,1) = \mathfrak{k} \oplus \mathfrak{p} $.

We choose our maximal abelian subalgebra of $\fp$ to be $\fa$ which is the linear span of $\left\{H= \left(
                                                                 \begin{array}{cc}
                                                                   0 & 1 \\
                                                                   1 & 0 \\
                                                                 \end{array}
                                                               \right)\right\}.$

Notice that for each $H \in \fa, \exp (tH) = \left(
                           \begin{array}{cc}
                             \cosh t & \sinh t \\
                             \sinh t & \cosh t \\
                           \end{array}
                         \right)$ for $t\in \R$, and that these matrices form the group $A$.

A straightforward calculation shows that the root spaces for $\fa$ are spanned by $X_\alpha = \left(
                                                                                                \begin{array}{cc}
                                                                                                  i & -i \\
                                                                                                  i & -i \\
                                                                                                \end{array}
                                                                                              \right)$
with root $\alpha(H)= 2$, and $X_{-\alpha} = \left(                                     \begin{array}{cc}
                                                                                                  -i & -i \\
                                                                                                  i & i \\
                                                                                                \end{array}
                                                                                              \right)$
with root $-\alpha$.

The weight $\rho$ is equal to $\frac{\alpha}{2}$, so that $\rho(H) =1$.

One calculates that $\theta X_\alpha = X_{-\alpha}$, that $X_\alpha + \theta X_\alpha = 2\left(
                                                                                           \begin{array}{cc}
                                                                                             i & 0 \\
                                                                                             0 & -i \\
                                                                                           \end{array}
                                                                                         \right)$ belongs to $\fk$, and that
 $X_\alpha - \theta X_\alpha =   2\left(
                                    \begin{array}{cc}
                                      0 & i \\
                                      -i & 0 \\
                                    \end{array}
                                  \right)$ belongs to $\fp$.

A simple calculation shows that $[X_\alpha, \theta X_\alpha] = 4H$.

 The Iwasawa decomposition for $SU(1,1)$ is $G=NAK$, where $K$ and $A$ are as above, and $N = \{ \exp tX_{\alpha}: t \in \R \}$ is the subgroup
$\left\{\left(
   \begin{array}{cc}
     1+it & -it \\
     it & 1-it \\
   \end{array}
 \right)  : t \in \R \right\}.$

 The subgroup $M$ of $K$ which fixes $H$ is $\{I, -I\}$. There are two irreducible representations of $M$, the trivial representation $\sigma_+(m) = 1$ for $m = \pm I$, and  the self representation $\sigma_-(\pm I) = \pm 1$.

 We now describe the principal series representations of $SU(1,1)$. As indicated in the appendix, they are indexed by the characters of $A$, specified by the real numbers $\lambda$, corresponding to the characters $\exp(tH) \mapsto e^{i\lambda t}$.

 For each $\lambda \in \R$, we define
 $\xi_{\lambda}$ acting on $L^2(SO(2))$ by

$$(\xi_{\lambda}(g)f)(l) = e^{(i\lambda + 1)(A(lg))}f(u(lg)),$$
where $f \in L^{2}(SO(2)), l \in SO(2), g = n\exp(A(g))u(g) \in SU(1,1)$.

The two subspaces of $L^2(SO(2))^{\pm}$ consisting (respectively) of the even functions $\{f: f(-k)=f(k)\}$ and the odd functions $\{f: f(-k)=-f(k)\}$ on the circle  then each carry an irreducible representation $\xi^+_{\lambda}$ and $\xi^-_{\lambda}$ , corresponding to the two irreducible representations of $M$.

Now we can decompose $\mathcal{H} = L^2(SO(2))$ into $K$-types, using the Fourier transform: that is, for all $n \in \mathbb{Z}$ we let $ \mathcal{H}_n$ be the linear span of the character $\theta \mapsto e^{in \theta}$; these spaces are all one-dimensional for $SU(1,1)$. We have $\mathcal{H}^{\pm} = \oplus_n \mathcal{H}_n$, the sum being over even integers for $\mathcal{H}^{+}$ and over odd integers for $\mathcal{H}^{-}$.

The generalised Eisenstein integrals are indexed by $\lambda$ together with a pair $(n',n)$ in $\Z \times \Z$.
\begin{equation}\label{gE}
\Phi_{\lambda,n', n}(g) = \frac{1}{2\pi}\int_0^{2\pi}e^{(i\lambda + 1)(A(k(\theta)g))}e^{in'(u(k(\theta)g))-in \theta} d\theta.
\end{equation}

The usual Eisenstein integrals are recovered when $n'=0$.

Now let us pass to hyperbolic space $\mathcal{D}$, which we may realise as $G/K=SU(1,1)/SO(2)$ via $M.z = \ds\frac{az+b}{\bar{b}z + \bar{a}}$, where of course
$M$ is the matrix $\left(                                                       \begin{array}{cc}
                                                                                   a & b \\
                                                                                   \bar{b} & \bar{a}\\
                                                                                 \end{array}
                                                                               \right) $.

Notice that the stabiliser of the point $z=0$ is precisely the group $K$. Thus the mapping $M \mapsto \frac{b}{\bar{a}}$ is a one-one correspondence between $G/K$ and $\mathcal{D}$.  Note that $\exp(tH)$ maps to $\mathrm{tanh}(t)$ through this correspondence. As shown by Helgason \cite{Hel2} (see section 4 of the Introduction, pp. 29 -72) the action of $SU(1,1)$ is isometric for the usual Riemannian structure of the Poincar\'e disc. The boundary $\mathcal{B}= \partial \mathcal{D}$ is identified as $K=SO(2)$.  Furthermore, consider the expression $A(k(\cdot)g)$ in the exponent of the integrand of formula (\ref{gE}). As the $A$ part of $g$ coincides with the $A$ part of $gk$ for all $g \in G,k \in K$, this expression passes to the quotient and defines an element of $A$ as a function of $k$ and $gK$. Writing this element of $A$ in the form $\exp(tH)$, let us denote $t=<gK, k>$. This is consistent with Helgason's notation $<z,b>$ for $z \in \mathcal{D}, \,\, b \in \mathcal{B}$.

Now we can express the image of the Eisenstein integrals (\ref{gE}) under the quotient map from $G$ to $\mathcal{D}$, as
\begin{equation*}
   \Phi_{\lambda, n}(z) := \Phi_{\lambda, 0, n}(z) = \int _{\mathcal B} e^{(i\lambda +1)<z,b>}b^{-n}db
\end{equation*}

It follows readily from this that if $z=e^{i\theta }r$, with $r\geq 0$, then $$\Phi_{\lambda, n}(z)= e^{i n\theta}\Phi_{\lambda, n}(r).$$

Furthermore, the following formula is derived in \cite{Hel2}, p.60:
\begin{equation}\label{HG}
    \Phi_{\lambda, n}(r) = r^{|n|}\frac{\Gamma(|n| + \nu)}{\Gamma(\nu)|n|!}{}_2F_1\left(\nu, 1-\nu;|n|+1;\frac{r^{2}}{r^{2}-1}\right),
\end{equation}
where $\nu = \frac{1}{2}(1 + i\lambda)$ and ${}_2F_1$ is the Gauss hypergeometric function.

Notice that $\Phi_{\lambda,0} (z) = \int_{\mathcal B} e^{(i\lambda +1)<z,b>}db$ is the usual spherical function, denoted $\phi_{\lambda}$ by Helgason.

Now if we take the Fourier transform of $f \in L^{1}(D)$, we obtain
\begin{eqnarray*}
  \int_{\mathcal{D}}f(z)\Phi_{\lambda, n}(z)dz &=& \int_{\mathcal{D}}f(z)\int _{\mathcal B} e^{(i\lambda +1)<z,b>}b^{-n}db dz\\
   &=&  \int _{\mathcal B}\int_{\mathcal{D}}f(z)e^{(i\lambda +1)<z,b>}dzb^{-n}db\\
   &=& \int _{\mathcal B}\tilde{f}(\lambda,b) b^{-n}db\\
   &=&  \widehat{\tilde{f}}(\lambda, \cdot)(n),
\end{eqnarray*}
where, in Helgason's notation, $\tilde{f}(\lambda,b)= \int_{\mathcal{D}}f(z)e^{(i\lambda +1)<z,b>}dz$ is the ``Fourier transform'' on $\mathcal{D}$.

Notice that Helgason's Fourier transform is a function of $(\lambda,k)$, for $k$ in the circle group. The Helgason Fourier transform can be recovered from the group Fourier transform by use of the usual Fourier series, specifically: $$\tilde{f}(\lambda,b)= \sum_{n \in \Z}\left(\int_{\mathcal{D}}f(z)\phi_{\lambda, n}(z)dz\right)b^n.$$
With this notation, the characteristic exponent in (\ref{LKMan}) has the form
$$\psi_{\lambda,n} = -b^i \rho_{\lambda,0, n}(X_i) + \sum_{k \in \Z} a^{i,j}\rho_{\lambda, 0, k}(X_i)\rho_{\lambda,k,n}(X_j) + \eta_{\lambda,n},$$
for each $\lambda  \in \R, n \in \mathbb{Z}$ and where we note that $i,j \in \{1,2\}$.

We can calculate $\rho_{\lambda, 0, k}$ as follows. The tangent space is identified with $\fp$, and we choose as orthonormal basis
$X_1 =  \frac{1}{\sqrt{2}}\left(\begin{array}{cc}                                          0 & 1 \\
                                                                   1 & 0 \\
                                                                 \end{array}
                                                               \right)$ and
                                                               $X_2  = \frac{1}{\sqrt{2}}\left(
                                                                         \begin{array}{cc}
                                                                           0 & i \\
                                                                           -i & 0 \\
                                                                         \end{array}
                                                                       \right)$.

Then $\exp tX_1 =\left(
                   \begin{array}{cc}
                     \cosh  \frac{t}{\sqrt{2}} & \sinh  \frac{t}{\sqrt{2}} \\
                     \sinh  \frac{t}{\sqrt{2}} & \cosh  \frac{t}{\sqrt{2}} \\
                   \end{array}
                 \right)$ which maps under the quotient map to $\tanh  \frac{t}{\sqrt{2}}$ in $\mathcal{D}$, and
$\exp tX_2 = \left(
               \begin{array}{cc}
                 \cosh  \frac{t}{\sqrt{2}} & i \sinh  \frac{t}{\sqrt{2}} \\
                 -i \sinh  \frac{t}{\sqrt{2}} & \cosh  \frac{t}{\sqrt{2}} \\
               \end{array}
             \right)$ which maps under the quotient map to $i \tanh \frac{t}{\sqrt{2}}$ in $\mathcal{D}$.
From the above, it follows that $\rho_{\lambda, 0, n} (X_2) = i^n \rho_{\lambda, 0, n}(X_1)$, and by definition,
  $\rho_{\lambda, 0, n}(X_1)= \frac{d}{dt}\Phi_{\lambda,n}(\tanh \frac{t}{\sqrt{2}}) |_{t=0}$. By equation (\ref{HG}),
  $$\Phi_{\lambda, n}(\tanh \frac{t}{\sqrt{2}})= (\tanh \frac{t}{\sqrt 2})^{|n|}\frac{\Gamma(|n| + \nu)}{\Gamma(\nu)|n|!}{}_2F_1\left(\nu, 1-\nu;|n|+1;\sinh^2 \frac{t}{\sqrt{2}}\right).$$
A simple computation involving the hypergeometric series shows that the derivative of this function evaluated at $t=0$ is $\frac{\nu(1-\nu)}{2} = \frac{1}{8}(1 + \lambda^2)$ if $n=0$,  $\frac{\nu}{2} = \frac{1}{4}( 1 + i \lambda)$ if $n = \pm 1$ and zero otherwise.

Thus $\rho_{\lambda, 0, 0}(X_1) = \rho_{\lambda, 0, 0}(X_2) = \frac{1}{4}(1 + \lambda^2) $,  $\rho_{\lambda, 0, n}(X_1) = \frac{1}{4}( 1 + i \lambda)$ if $n=\pm 1$, and  $\rho_{\lambda, 0, n}(X_2) = \frac{\pm i}{4}( 1 + i \lambda)$ if $n=\pm 1$. Both coefficients are zero for $n \neq \pm 1,0$.

These formulae are the expressions for hyperbolic space of the infinitesimal action of Lemma \ref{diffrep}.

We have  $\eta_{\lambda,0} = \int_{\mathcal{D}}(\phi(z)-1)\tilde{\nu}(dz)$, and for $n \neq 0$
\begin{equation*}
   \eta_{\lambda,n} = \int_{\mathcal{D}} (\Phi_{\lambda,n}(z) + \tilde{x_i}(z)\rho_{\lambda, 0, n}(X^{i}))\tilde{\nu}(dz),
\end{equation*}
where $\tilde{\nu}$ is a L\'{e}vy measure on $D$ defined in terms of the given L\'{e}vy measure $\nu$ on $SU(1,1)$ by
$\tilde{\nu} = \nu \circ \natural^{-1}$. Similarly $\tilde{x_i}(z) = x_{i} \circ \natural^{-1}$ for $i=1,2$.
In the spherically symmetric case, the spherical function $\phi$ is a Legendre function and an explicit L\'{e}vy-Khintchine formula is presented in Theorem 7.4 of \cite{Get}. To recover this result, notice that for the $K$-bi-invariant case, we have \bean \psi_{\lambda, 0} & = & a \Sigma_{i=1,2} \rho_{\lambda, 0, 0}(X_i)^2 + \eta_{\lambda,0}\\ & = & \frac{1}{2}a(1 + \lambda^2) + \int_{\mathcal{D}}(\phi(z)-1)\tilde{\nu}(dz) \eean
for a non-negative constant $a$.

\section{Appendix: Some Aspects of Representation Theory}\label{sec8}
In this paper, we have used the Fourier transform on a semisimple group as an essential tool. In order to make the article more self-contained and accessible, we give here a brief introduction to this subject, showing how the Mackey theory of induced representations (which derives from work of Schur) is applied to the case of semisimple Lie groups to derive the principal series, and giving some of the basic ideas about the Fourier transform for locally compact groups.
A certain divergence of notation can be found in the literature: one can choose whether to write the induced representation on the left (\`a la Wallach \cite{Wa}) or the right (\`a la Knapp \cite{Kn}); we have chosen the latter, which means there is no $g^{-1}$ in the formula, but means that the representation is realised on $L^2$ of the right coset space $H\backslash G$. One must also decide whether to write the Iwasawa decomposition as $G=KAN$ or $ANK$ or $NAK$. We have chosen $NAK$. Each of the other choices gives slightly different formulae, although of course the representation theory is the same. Here we have aimed at a consistent choice that provides simple formulae, consistent with those of Helgason \cite{Hel3}. Notice that although we use right coset spaces for the representation theory of $G$, our choice of Iwasawa decomposition means that it is natural to represent the symmetric space $G/K$ as \textbf{left} cosets, as we do in section \ref{sec5} above.

\subsection{Induced Representations}
Let $G$ be a locally compact group and $H$ a closed subgroup. Mackey's theory shows how to start from a (continuous unitary) representation $\pi$ of $H$ acting in $\mathcal{H}_\pi$ and construct a continuous unitary representation $\xi = \pi \uparrow_H^G$ of $G$, which generalises Schur's construction from finite groups, and has good functorial properties.

Before making the construction, we state Mackey's Lemma.
\begin{lemma} \label{Mack}If $G$ is second countable and $H$ is closed in $G$, then the quotient space $H\backslash G = \{Hg: g \in G\}$ has a right action of $G$ by $g : Hg_1 \mapsto Hg_1g^{-1}$ and there exists a Borel measure $\mu$ on $H\backslash G$ which is quasi-invariant for this action, with the property that the measure of each open set is strictly positive. Furthermore there is at least one open set having finite $\mu$-measure.
\end{lemma}
By Lemma \ref{Mack}  and the Radon-Nikod\'ym Theorem, we can define an a.e. unique function $\frac{d\mu \circ g}{d\mu}: H\backslash G \to \R^+$ such that for all Borel sets $E \subseteq H\backslash G$, $$\int_E \frac{d\mu \circ g}{d\mu}(h) \mu(dh) = \mu(Eg^{-1}).$$

The construction will give $\xi$ acting by a slight modification of the right regular representation in a space of functions on $G$.

 Given $\pi$ and $\mathcal{H}_\pi$ as above, we let
\begin{equation}\label{indref}
    \mathcal{H}_{\xi} = \{ f: G \to \mathcal{H}_\pi : f(hg) = \pi(h)f(g) \,\, \mathrm{and} \,\,  \int_{H\backslash G} \parallel f(Hg) \parallel^2 d\mu(Hg) < \infty \}
\end{equation}
In order to make sense of the integral in the above definition, notice that if $f$ and $k$ both satisfy the invariance condition, then their inner product is constant on cosets of $H$ in the sense that
\begin{equation}
    <f(hg), k(hg)>_{\mathcal{H}_\pi} = <f(g), k(g)>_{\mathcal{H}_\pi}.
\end{equation}
Thus, taking $f=k$, the integrand in (\ref{indref}), evaluated at $g$, is seen to depend only on the coset $Hg$. In fact, we may define an inner product on ${\mathcal{H}_\xi}$ by
$<f, k>_{\mathcal{H}_\xi} = \int_{H\backslash G} <f(Hg), k(Hg)> d\mu(Hg).$

Now, we define $\xi(g)f(g_1) = f(g_1g) \left(\frac{d\mu\circ g}{d\mu}(Hg_1)\right)^{1/2}$.  A simple calculation shows that $\xi$ is unitary.

Note that each function in $\mathcal{H}_{\xi}$ is determined by its values on a section of the quotient $H\backslash G$. Indeed, given a section $\gamma: H\backslash G \to G$, i.e. a measurable map such that $\gamma(Hg) $ lies in the same coset as $g$, we can uniquely factorise each element $g$ of $G$ as $g= h(g) \gamma (Hg)$ for some $h(g) \in H$. Now if $f$ belongs to $\mathcal{H}_{\xi}$, then $$f(g)= f(h(g) \gamma(Hg) = \pi(h(g) f(\gamma(Hg)).$$

Actually, it is better in some cases to think of $\mathcal{H}_{\xi}$ as the space $L^2(H\backslash G, d\mu) \bigotimes \mathcal{H}_\pi$. In this picture, the representation $\xi$ acts on an elementary tensor
$f\otimes \eta$, $f \in L^2(H\backslash G, d\mu)$ and $\eta \in \mathcal{H}_\pi$ by
\begin{equation}\label{sectpic}
    \xi(g)f\otimes \eta = \mathcal{R}_g f \otimes \pi(h(g)) \eta,
\end{equation}
where $\mathcal{R}_g$ is the natural action of $G$ on $L^2(H\backslash G)$ given by $\rho_gf(Hg_1) = f(g_1g)(\frac{d\mu\circ g}{d\mu}(Hg_1))^{1/2}$. The presence of the Radon-Nikod\'ym derivative ensures that the $L^2$ norm is preserved.

\subsection{The Principal Series}
We are going to apply the theory of the previous subsection to the case where $G$ is a semi-simple Lie group and $H$ is a minimal parabolic subgroup of $G$.

We have the Iwasawa decomposition $G=NAK$ and we write $g \in G$ as $g= n(g) \exp (H(g)) \tau(g)$, where $n(g) \in N$, $H(g) \in \mathfrak{a}$ and $\tau(g) \in K $.

As usual, $M$ is the centraliser of $A$ in $K$. The minimal parabolic subgroup of $G$ associated to these choices is $NAM$. As $M$ normalises both $A$ and $N$, this group has the structure of a semi-direct product of $N$ by the actions of $A$ and $M$. Given $\lambda \in \mathfrak{a}^*$, we have a character of $A$ given by $\exp(H) \mapsto e^{i \lambda(H)}$ for $H \in \mathfrak{a}$. This gives a character of $NA$ by $(n,\exp(H))\mapsto e^{i \lambda(H)}$. Given also an irreducible representation $\nu$ of $M$, we let $\pi_{\lambda,\nu}(n, \exp(H),m) \eta = e^{i \lambda(H)}\nu(m) \eta$, and this provides us with a representation of $NAM$ in $\mathcal{H}_\nu$.

The principal series is defined to be the family $\xi_{\lambda,\nu} = \pi_{\lambda,\nu}\uparrow _{NAM}^G$.

We give some more explicit formulae for the principal series representations. Firstly, notice that the coset space $H\backslash G$ is $NAM \backslash NAK \sim M\backslash K$. If we realise the coset space as $M\backslash K$, we have
\begin{lemma} \label{pf} The push-forward of the Haar measure of $K$ gives an $K$ invariant measure on $M\backslash K$. Furthermore, this measure is quasi-invariant for the action of $NAK $. The Radon-Nikod\'ym derivative for the action of $G$ is given by $$\frac{d\mu \circ g}{d\mu}(Mk) = e^{2\rho(H(kg))},$$where $\rho$ is half the sum of the positive roots of $(\mathfrak{g}, \mathfrak{a})$.
\end{lemma}

\vspace{5pt}

Now, using the formulae of the previous section, $\xi_{\lambda,\nu}$ acts in the space $$\mathcal{H}_{\lambda,\nu} = \{f: G=NAK \to \mathcal{H}_\nu : f(namg)= e^{i\lambda(H)}\nu(m)f(g)\}$$with the action given by $$\xi_{\lambda,\nu}(g)f (g_1) = f(g_1g) e^{\rho(H(g_1g))}.$$

It is interesting to re-write this action using the notation of formula (\ref{sectpic}), as we then obtain an explicit representation on $L^2(M\backslash K)$.
Notice that the map $g\mapsto \tau(g)$ maps $nak$ to $k$, and, because of the normalising property of $M$, $\tau(mg)= m\tau(g)$. It follows that $g\mapsto M\tau(g)$ is a section for the quotient map $g \mapsto NAMg$.

We may now write the action of $\xi_{\lambda, \nu}(g)$ on $L^2(M\backslash K) \bigotimes \mathcal{H}_\nu$ as $$\xi_{\lambda, \nu}(g)(f \otimes \eta)(l) = e^{(i\lambda + \rho)(H(lg))}\mathcal{R}_gf\otimes \nu(m(k))\eta. $$

Here, $\mathcal{R}_g f (Mk) = f(M \tau(kg))$, and we suppose that there is a section $p:M\backslash K \to K$ where each $k \in K$ can be written $k= p(Mk)m(k)$ for a suitable $m(k) \in M$.

Actually, there is a slightly more straightforward way of writing this representation. Consider the representation of $G$ on $L^2(K) $ given by
\begin{equation}\label{kversion}
    \xi'_{\lambda}(g)(f)(l) = e^{(i\lambda + \rho)(H(lg))}f(\tau(lg)),
\end{equation} for $g \in G, l \in K$.
This representation extends to $ \xi'_{\lambda}\otimes \textbf{1}$ acting on $L^2(K) \bigotimes \mathcal{H}_\nu$ trivially on the second coordinate. Now consider the subspace generated by the elementary tensors such that $\xi'_{\lambda}(m)f\otimes \eta = f\otimes \nu(m)\eta$ for all $m \in M$. The restriction of $ \xi'_{\lambda}\otimes \textbf{1}$ to this subspace is equivalent to $\xi_{\lambda, \nu}$.

Notice that $\xi'_{\lambda}(k)f(k_1) = f(k_1k)$ is exactly the right regular representation of $K$ on $L^2(K)$.

In fact, see \cite{Kn}, for each cuspidal parabolic subgroup $Q$ of $G$, there is a $Q$-series of representations of $G$; we have only dealt with the case when $Q$ is a minimal parabolic subgroup. Another feature that one can vary is the assumption that $\lambda$ is real-valued; thus we obtain non-unitary principal series. The {\it Casselman subrepresentation theorem} states that every irreducible admissible representation of $G$ appears (on the $K$-finite level) as a subrepresentation of some non-unitary principal series. The description of all irreducible representations of $G$ is thus somewhat more complicated than we have described here.

\subsection{K-types}

Suppose that $G$ is a locally compact group and $K$ a compact subgroup. We are going to apply this theory to the case where $G$ is a connected semi-simple Lie group having a finite centre and $K$ is a maximal compact group, but for the moment, we keep the discussion general.

If $\xi $ is a representation of $G$ on $\mathcal{H}$, we want to consider the restriction $\xi |_K$. This is a representation of a compact group, which may be written as a direct sum of irreducible representations of $K$. Letting $\hat{K}$ be the unitary dual of $K$ (i.e. a maximal set of pairwise inequivalent irreducible representations of $K$), we may write
\begin{equation}\label{mult}
\xi |_K = \bigoplus_{\pi \in \hat{K}} m_{\pi} \pi,
\end{equation}
where $1 \leq m_{\pi} \leq \infty$ is the multiplicity of $\pi$ in the decomposition.

If $G$ is a connected semi-simple Lie group having a finite centre, $K$ is a maximal compact subgroup and $\xi'$ is the principal series constructed in the previous section, we saw that $\xi'|_K$ is the regular representation of $K$ on $L^2(K)$. By the Peter-Weyl theorem, the regular representation contains each irreducible representation $\pi$ of $K$ exactly $d_{\pi}$ times. Thus, in this example, $m_{\pi}$ is nothing but the dimension of the representation space.

We shall continue our general discussion, under the assumption that $m_{\pi} < \infty$ for all $\pi \in \hat{K}$.
Equation (\ref{mult}) means that there exists a family of orthogonal projections on $\mathcal{H}$, $\{P^{\pi}_i :\pi \in  \hat{K}, i=1, \dots, m_{\pi}\}$, each with finite range of dimension $d_{\pi}$, so that $P^{\pi}_i \xi|_K P^{\pi}_i$ is equivalent to $\pi \in \hat{K}$, for $i=1,\dots, m_{\pi}$, and such that $$\bigoplus_{\pi \in \hat{K}} \bigoplus_{i=1}^{m_{\pi}} P^{\pi}_i = \mbox{Id}.$$

Let $P^{\pi}_i\mathcal{H} = \mathcal{H}^{\pi}_i$. Then $\mathcal{H}^{\pi}_i$ carries an irreducible representation of $K$ equivalent to $\pi$ for $i=1,\dots, m_{\pi}.$ Furthermore, this representation acts by $\pi(k)\phi = \xi(k) \phi$ for $\phi \in \mathcal{H}^{\pi}_i$.

This decomposition of $\mathcal{H}$ is called a decomposition into \textbf{$K$-types}.

\subsection{Fourier Transforms and Convolutions}\label{sec74}

Now suppose that $\mu$ is a Borel probability measure on a locally compact group $G$, and $\xi$ is an irreducible representation of $G$ acting in $\mathcal{H}_{\xi}$. Then we may define the Fourier transform of $\mu$ at $\xi$ by:
$$\hat{\mu}(\xi) = \int_G \xi(g^{-1}) \mu(dg).$$

This is to be thought of as a Bochner integral and defines an operator in $\mathcal{L}(\mathcal{H}_{\xi})$ by $$\langle \hat{\mu}(\xi)v, w\rangle = \int_G \langle \xi(g^{-1})v,w\rangle \mu(dg),$$ for $v,w \in \mathcal{H}_{\xi}$.

In section \ref{sec3} we use generalised Eisenstein integrals to express this operator in terms of its projections on the $K$-types.

If $\mu_{1}$ and $\mu_{2}$ are Borel probability measures on $G$, their convolution $\mu_{1} * \mu_{2}$ is a Borel probability measure defined on $G$ for which
\begin{equation} \label{conv1}
\int_{G}f(g)(\mu_{1} * \mu_{2})(dg) = \int_{G}\int_{G}f(gh) \mu_{1}(dg)\mu_{2}(dh),
\end{equation}
where $f$ is an arbitrary bounded Borel measurable function on $G$.

It is now standard and easy to see that
$$ (\widehat{\mu^{(1)} * \mu^{(2)}})({\xi}) =   \widehat{\mu^{(1)}(\xi)}\widehat{\mu^{(2)}(\xi)},$$
where the right hand side is composition of operators in $\mathcal{L}(\mathcal{H}_{\xi})$.

\vspace{5pt}

{\it Acknowledgements}

\vspace{5pt}
We acknowledge the support of the Australian Research Council. DA would like to thank the University of New South Wales for generous hospitality and support during January-March 2012 when much work was carried out on this paper. Both authors are grateful to Ming Liao for very helpful comments on an earlier version of this paper. We would also like to thank the referees for highly valuable feedback which has greatly improved the manuscript.

\end{document}